\documentclass[12pt]{article}

\usepackage[T2A]{fontenc}

\usepackage{amsfonts}
\usepackage{amssymb}
\usepackage{amsmath}
\usepackage{amsthm}

\def \le {\leqslant}
\def \ge {\geqslant}

\begin{document}

\centerline{\bf Differentiability of the Minkowski question mark function.} \vskip+0.5cm \centerline{\bf  Anna A. Dushistova, Igor D. Kan,
Nikolai G. Moshchevitin \footnote{ Research is supported by   RFBR  grant No. 09-01-00371-a
 } } \vskip+0.5cm

\begin{abstract}
 We prove   new results on the derivative of the Minkowski question mark function.
Some of our theorems are  non-improvable.
 \end{abstract}

{\bf 1.\,\,\, The Minkowski function  $?(x)$.} \,\,\, The function $?(x)$ is defined as follows.
 $ ?(0)=0, ?(1)= 1$,
 if the values $  ?\left(\frac{p}{q}\right)$ and $?\left(\frac{p'}{q'}\right)$ are defined for consecutive  Farey fractions
$\frac{p}{q}, \frac{p'}{q'}$ then
$$?\left(\frac{p+p'}{q+q'}\right) = \frac{1}{2} \left( ?\left(\frac{p}{q}\right)+
?\left(\frac{p'}{q'}\right)\right);
$$
for irrational $x\in [0,1]$ function $?(x)$ is defined by continuous arguments. This function firstly was considered by H. Minkowski  (see. \cite{MI},
p.p. 50-51) in 1904. $?(x)$ is a continuous increasing function. It has derivative almost everywhere. It satisfies Lipschitz condition
\cite{SA}, \cite{KI}. It is a well-known fact that the derivative  $?'(x)$ can take only two values -  $0$  or $+\infty$. Almost everywhere we
have
 $?'(x)=0$.
Also if irrational $x=[0;a_1,...,a_t,...] $ is represented as a regular  continued fraction with natural partial quotients then
$$
?(x) = \frac{1}{2^{a_1-1}} - \frac{1}{2^{a_1+a_2-1}}+ ...+ \frac{(-1)^{n+1}}{2^{a_1+...+a_n-1}}+ ...  .
$$
These and some other results one can find for example in papers \cite{SA},\cite{PARA},\cite{PARA2}.
Here we should note the connection between function $?(x)$ and Stern-Brocot sequences. We remind the reader the definition of Stern-Brocot
sequences $ F_n $, $ n=0, 1, 2, \dots $. First of all let us put
  $ F_0=\{0, 1\}=\{\frac{0}{1}, \frac{1}{1}\} $.
Then for the sequence $ F_n $ treated as increasing sequence of rationals $
 0=x_{0, n}<x_{1, n}< \dots <x_{N \left( n
  \right), n}=1, N(n)=2^{n},
$  $x_{j, n} =p_{j, n}/q_{j, n}, ( p_{j, n},q_{j, n}) = 1$ we
define the next sequence $F_{n+1}$ as  $ F_{n+1} = F_n \cup
Q_{n+1} $ where  $ Q_{n+1} $ is the set of the form $ Q_{n+1}=\{
x_{j-1,n}\oplus x_{j,n}, j=1, \dots , N(n)\}.  $ Here operation
$\oplus$ means taking the mediant fraction for  two rational
fractions:
  $ \frac{a}{b}\oplus \frac{c}{d} = \frac{a+c}{b+d}$.
The Minkowski question mark function $?(x)$ is the limit distribution function for Stern-Brocot sequences:
$$
?(x)  =\lim_{n\to \infty} \frac{\# \{ \xi\in F_n :\,\,\,\xi \le x\}}{2^n +1}.
$$

 {\bf 2.\,\,\, Notation and parameters.} \,\,\,
 For natural numbers $a_1,...,a_t$ the notation $\langle a_1,...,a_t \rangle $ denotes the continuant with digits $a_1,...,a_t$.
 That is
empty continuant is equal to one,
$$
\langle a_1\rangle=a_1,\,\, \langle a_1,...,a_t \rangle = a_t \cdot \langle a_1,...,a_{t-1}\rangle+ \langle a_1,...,a_{t-2}\rangle ,\,\, t \ge
2.$$

For irrational  $x \in (0,1)$ we consider the continued fraction expansion
$$
x= [0;a_1,a_2,...,a_t,...]= \frac{1}{\displaystyle{a_1+\frac{1}{\displaystyle{a_2 + \dots+\frac{1}{\displaystyle{a_t + \displaystyle{\dots} }}}}}}
$$ with natural partial quotients $a_t$.
For breavity we use the notation
$$ x =[a_1,...,a_t,...].$$

  For
 the continued fraction under consideration the convergent fraction of order  $t$ is denoted as $p_t/q_t= p_t(x)/q_t(x) = [a_1,a_2,...,a_t]$ (hence,
 $ q_t= q_t(x)
=\langle a_1,a_2, ...,a_t\rangle,   p_t =p_t(x)=\langle a_2,...,a_t\rangle $). For an irrational number $x$ we consider the sum $S_x(t)$ of its
partial quotients up to $t$-th:
$$
S_x(t) = a_1+a_2+...+a_t.
$$

 We need numbers
\begin{equation}
\lambda_j = \frac{j+\sqrt{j^2+4}}{2} ,\,\,\, j =1,2,3,..., \label{lambdaj}
\end{equation}
\begin{equation}
\mu_j = \frac{j+2+\sqrt{j^2+4j}}{2} ,\,\,\, j =1,2,3,..., \label{muj}
\end{equation}
\begin{equation}
\kappa_1 = \frac{2\log \lambda_1}{\log 2}= \frac{2\log  \frac{1+\sqrt{5}}{2}}{\log 2} = 1.388^+, \label{kappa1}
\end{equation}
\begin{equation}
\kappa_2 = \frac{4\log\lambda_5-5\log\lambda_4}{\log \lambda_5-\log\lambda_4-\log\sqrt{2}}= \frac{4\log
\frac{5+\sqrt{29}}{2}-5\log(2+\sqrt{5})}{\log \frac{5+\sqrt{29}}{2}-\log(2+\sqrt{5})-\log\sqrt{2}}= 4.401^+. \label{kappa2}
\end{equation}
Also for $n \ge 5$  we need the number
\begin{equation}
\kappa_{ n} = \frac{(n+1)\log \frac{1+\sqrt{5}}{2}-\log \frac{n+2+\sqrt{n^2+4n}}{2}}{ (n-1) \log \sqrt{2} - \log \frac{n+2+\sqrt{n^2+4n}}{2} + 2
\log \frac{1+\sqrt{5}}{2}} =$$ $$ = \frac{(n+1)\log\lambda_1-\log \mu_n}{ (n-1) \log \sqrt{2} - \log \mu_n + 2 \log \lambda_1}
 =\kappa_1+
\frac{4\log \lambda_1-2\log 2}{(\log 2)^2}\cdot \frac{\log n}{n}
+O\left(\left(\frac{\log n}{n}\right)^2\right) ,\,\, n \to \infty . \label{kappa1n}
\end{equation}

{\bf 3.\,\,\, A result by J. Paradis,  P. Viader, L. Bibiloni.}\,\,\, In \cite{PARA2} the following statement is proved.

{\bf Theorem A.}\,\,\,{\it

1. Let for real irrational $x\in (0,1)$   with $\kappa_1 $ from (\ref{kappa1}) one has
$$
 \limsup_{t\to
\infty}\frac{S_x(t)}{t} <\kappa_1 .$$ Then if
 $?'(x)$ exists the equality $?'(x)=+\infty $  holds.

 2. Let $\kappa_3 = 5.319^+$
be the root of equation
 $ 2\log (1+z )=z\log 2 $.
 Let for real irrational $x\in (0,1)$ one has
$$
 \liminf_{t\to
\infty}\frac{S_x(t)}{t} \ge \kappa_3.$$ Then if
 $?'(x)$ exists the equality $?'(x)=0$  holds.
 }

{\bf 4.\,\,\, Our main  results.} \,\,\, In this paper we prove the following theorems. In Theorems 1 -- 4 below $\kappa_1$ is taken from
(\ref{kappa1}) and $\kappa_2$ is taken from (\ref{kappa2})

{\bf Theorem 1.}\,\,\,{\it

(i) Let for an irrational number $x$ there exists a constant $C$
such that for all natural $t$ one has
\begin{equation}
S_x(t)\le \kappa_1t+\frac{\log t}{\log 2}+C. \label{eq1}
\end{equation}
Then $?'(x)$ exists and $?'(x)=+\infty$.

(ii) Let $\psi (t)$ be an increasing function such that $\lim_{t\to +\infty} \psi (t) =+\infty$. Then there exists an  irrational number $x\in
(0,1)$ such that $?'(x)$ does not exist and for any $t$ one has
\begin{equation}
S_x(t)\le \kappa_1t+\frac{\log t}{\log 2} +\psi(t). \label{eq2}
\end{equation}}

{\bf Theorem 2.}\,\,\,{\it

(i) Let for an irrational number $x\in (0,1)$ the derivative
$?'(x)$ exists and $?'(x)=0$. Then
for any real function $\psi =\psi (t)$ under conditions
$$
\psi (t) \ge 0,\,\,\,\, \psi (t) = o\left(\frac{\log\log t}{\log t}\right),\,\, t\to\infty
$$
there exists $T $ depending on $\psi$ such that for all $t\ge T$ one has
 $$
\max_{u\le t}\left(S_x(u)-\kappa_1 u\right) \ge \frac{ \sqrt{2\log \lambda_1-\log 2} }{\log 2}\cdot\sqrt{t\log
t}\cdot (1-t^{-\psi (t)}).
$$

(ii) There exists an irrational  $x\in (0,1)$  such that $?'(x)=0$ and for all  $t$ large enough one has
$$
S_x(t)-\kappa_1 t \le
 \frac{ \sqrt{16\log \lambda_1-8\log 2} }{\log 2}\cdot\sqrt{t\log
t}\cdot \left(1+2^5\left(\frac{\log\log t}{\log t}\right)\right).
$$

}

{\bf Theorem 3.}\,\,\,{\it

(i) Let for an irrational number $x$ there exists a   constant $C
$ such that for all natural $t$   one has
\begin{equation}
S_x(t)\ge \kappa_2t-C. \label{eq3}
\end{equation}
Then $?'(x)$ exists and $?'(x)=0$.

(ii) Let $\psi (t)$ be an increasing function such that $\lim_{t\to +\infty} \psi (t) =+\infty$. Then there exists an irrational number $x\in
(0,1)$ such that $?'(x)$ does not exist and and for any $t$ one has
\begin{equation}
S_x(t)\ge \kappa_2t-\psi(t). \label{eq4}
\end{equation}}

{\bf Theorem 4.}\,\,\,{\it

(i) Let for an irrational number $x\in (0,1)$ the derivative $?'(x)$ exists and $?'(x)=+\infty$. Then for  any  $t$ large enough one has

\begin{equation}
\max_{u\le t}\left( \kappa_2 u -S_x(u)\right) \ge  \frac{ \sqrt{t}}{10^8}.
\label{theoche}
\end{equation}

(ii) There exists an irrational  $x\in (0,1)$  such that $?'(x)=+\infty$ and for   $t$ large enough one has
$$
\kappa_2 t - S_x(t)\le 200 \sqrt{t}.
$$
}

From Theorems 1 -- 4 we immediately deduce the following result
which is stronger than the result announced by the first and the
third authors in the preprint \cite{arxiv}. We should note that
the statements of the next Corollary (and hence Theorems 1 -- 4  )
improves  Theorem A  by J. Paradis, P. Viader and L. Bibiloni
cited in Section 3.

 {\bf Corollary 1.}{\it

1. Let for real irrational $x\in (0,1)$   one has
$$
 \limsup_{t\to
\infty}\frac{S_x(t)}{t} <\kappa_1 .$$ Then
  $?'(x)$ exists and $?'(x)=+\infty $.

2. There exists an  irrational
 $x$ such that
$$\lim_{t\to \infty}\frac{S_x(t)}{t}=\kappa_1 $$
and
 $?'(x)=0$.

3. Let for real irrational $x\in (0,1)$  one has
$$
 \liminf_{t\to
\infty}\frac{S_x(t)}{t} >\kappa_2. \label{liminf}
$$
Then
  $?'(x)$ exists and $?'(x)=0$.

4. There exists an   irrational
 $x$ such that
$$\lim_{t\to \infty}\frac{S_x(t)}{t} = \kappa_2  $$
and
 $?'(x)=+\infty$.}

 {\bf Remark.}\,\,\, It is possible to prove that for any $\lambda$ from the interval
 $$
 \kappa_1\le \lambda\le \kappa_2
 $$
 there exist irrationals $x,y,z \in [0,1]$ such that
 $$
\lim_{t\to \infty}\frac{S_x(t)}{t}=\lim_{t\to \infty}\frac{S_y(t)}{t}=\lim_{t\to \infty}\frac{S_z(t)}{t}=\lambda
$$
and $?'(x) =0, ?'(y) = +\infty$ but $?'(z)$ does not exist.

 {\bf 5.\,\,\, Results on numbers with bounded partial quotients.} By $E_n, n \ge 2$ we denote the   set of
irrational numbers $x$  from the interval $(0,1)$ such that in the continued fraction expansion $ x= [a_1,a_2,...,a_t,...]$
 all partial quotients $a_j$ are bounded by $n$:
 $$
 a_j \le n ,\,\,\ \forall  j=1,2,3,... .
$$
In this paper we prove the following three theorems about the values of the derivative of the function $?(x)$  when $x \in E_n$. In theorems 5,6
below $\kappa_{n}$ is taken from (\ref{kappa1n}).

 {\bf Theorem 5.}\,\,\,{\it

(i) Let $ n \ge 5$ and $ x \in E_n$. Let for some constant $ C $ and for any natural $t$ one has
\begin{equation}
S_x(t)\le \kappa_{n} t+C. \label{eq1n}
\end{equation}
Then $?'(x)$ exists and $?'(x)=+\infty$.

(ii) Let $\psi (t)$ be an increasing function such that
$\lim_{t\to +\infty} \psi (t) =+\infty$. Then for a given $n\ge 5$
there exists an irrational number $x\in E_n$ such that $?'(x)$ does
not exist and for any natural $t$ one has
\begin{equation}
S_x(t)\le \kappa_{n}t+  \psi(t). \label{eq2n}
\end{equation}}

{\bf Theorem 6.}\,\,\,{\it

(i) Let for   $x \in E_n, n\ge 5$ the derivative $?'(x)$ exists and $?'(x)=0$. Then for   $t$ large enough one has
$$
\max_{u\le t}\left(S_x(u)-\kappa_{n}u\right) \ge  \frac{\sqrt{t }}{160 (n+2)^{10}}.
$$

(ii) For a given $n\ge 5$ there exists   $x\in E_n$  such that $?'(x)=0$ and for  $t$ large enough one has
$$
S_x(t)-\kappa_{ n}t \le 15n\sqrt{t }.
$$
}

{\bf Theorem 7.}\,\,\,{\it Let     $x \in E_4 $. Then $?'(x)$
exists and $?'(x)=+\infty$.}

{\bf 6.\,\,\, Lemma about the derivative of the function $?(x)$.} \,\,\, First of all we should note that the set $F_n$ consists of all
fractions $\frac{p}{q}\in [0,1]$ such that in the continued fraction expansion $ \frac{p}{q} =[a_1,...,a_t]$ one has $a_1+...+a_t \le n$.  As
for the set $Q_n =\ F_n\setminus F_{n -1}$ it consists of all fractions $\frac{p}{q}\in [0,1]$ satisfying the condition $a_1+...+a_t = n$.

For a real $x =[a_1,...,a_r,...   ]$ we consider two convergents $
\frac{p_{r-1}}{q_{r-1}} =[a_1,...,a_{r-1}]$ and $ \frac{p_r}{q_r}
=[a_1,...,a_r]$. As $\left|
\frac{p_{r-1}}{q_{r-1}}-\frac{p_{r}}{q_{r}}\right| =
\frac{1}{q_rq_{r-1}} $ and fractions are consecutive elements from
$F_n$ (where $n =a_1+...+a_r$) we see that
\begin{equation}
   \frac{\left|?\left( \frac{p_{r-1}}{q_{r-1}}\right) -?\left( \frac{p_{r}}{q_{r}}\right)\right|}{\left|
   \frac{p_{r-1}}{q_{r-1}}-\frac{p_{r}}{q_{r}}\right| }=
\frac{q_rq_{r-1}}{2^{a_1+...+a_{r} }}.
  \label{12}
\end{equation}
In this section we generalize the  equality (\ref{12}).

 We must do another useful observation. Let $\xi^0,\xi^1$ be  two consecutive elements of $F_n$ and $ x= [a_1,...,a_t,...] \in
(\xi^0,\xi^1)$. Consider the fraction $\xi = \xi^0\oplus \xi^1 \in F_{n+1}$ (notation $\oplus$ is defined in Section 1) and suppose $ x\neq
\xi$. In this situation one of the two fractions $ \xi^0,\xi^1$ lie on the same side from $x$  with $\xi$. Then the fraction on the opposite
side must be a convergent fraction to $x$.

 {\bf Lemma 1.}\,\,\,{\it For an
irrational $x= [a_1,...,a_t,...]\in (0,1)$ and    $\delta  $  small enough there exists  a natural
  $ r = r (x,\delta )$  such that
\begin{equation}
   \frac{?(x+\delta ) - ?(x)}{\delta }\ge
\frac{q_rq_{r-1}}{2^{a_1+...+a_{r}+4}}
  \label{DDDDw1}
\end{equation}
}

{\bf Proof.}\,\, It is enough to prove Lemma 1 for positive  $\delta$\label{SX1}.  Define the unique natural number  $n$ such that  $F_n \cap (x,x+\delta )
= \emptyset$, $F_{n+1} \cap (x,x+\delta ) = \xi$.
 So $(x,x+\delta )\subset (\xi^{0},\xi^1)$, where $ \xi^0 =\frac{p^0}{q^0}, \xi^1=\frac{p^1}{q^1}$ are two
 successive points from the finite set $F_n $. Then $\xi =\xi^0\oplus \xi^1 =\frac{p^0+p^1}{q^0+q^1}$.
 We see that $\xi$ and $\xi^1$  both lie on the same side from $x$. Then as it was mentioned above one can
 easily see from the Farey tree construction that for some natural $t$ will happen $\xi^0 = p_t(x)/q_t(x)=
p_t/q_t$.   At the same time rationals $\xi$ and $\xi^1$ must be among convergent fractions to $x$ or intermediate fractions to $x$
(intermediate fraction is a fraction of the form $\frac{p_ta+p_{t-1}}{q_t a+q_{t-1}}, 1\le a< a_{t+1}$).

Define natural $z$ to be minimal such that either  $ \xi_- = \xi^0\underbrace{\oplus \xi \oplus ... \oplus \xi}_{z} \in (x,\xi) $ or
 $\xi_{+} = \xi^1 \underbrace{\oplus \xi \oplus ... \oplus \xi}_{z} \in (\xi  , x+\delta)$.
 Then
$ \xi_{- -} = \xi^0\underbrace{\oplus \xi \oplus ... \oplus \xi}_{z-1} \le x $ and
 $\xi_{++} = \xi^1 \underbrace{\oplus \xi \oplus ... \oplus \xi}_{z-1} \ge  x+\delta$.
 As points $\xi_{--}<\xi_-<\xi <\xi_+<\xi_{++}$
are successive points from $F_{n+z+1}$
 and $?(x)$ increases, we have
\begin{equation}
\frac{1}{2^{n+z+1}} \le \min \{ ?(\xi_+) - ?(\xi ), ?(\xi ) - ?(\xi_-)\} \le
 ?(x+\delta ) - ?(x).
  \label{w2w1}
\end{equation}

As $\delta $ is small we can suppose that $n>1$. So  $q^0\neq q^1$. Then there may be two opportunities: $q^0>q^1$ ({\bf case 1}) and $q^0<q^1$
({\bf case 2}).

{\bf Case 1.} In this case we have $q_t = q^0>q^1$. We see that $ \xi^0 =\frac{p_t}{q_t} \in F_n\setminus F_{ n - 1}$ but $\xi^1 \in F_{ n -
1}$. So $ n =a_1+...+a_t$.
 Then
for some  fraction $\frac{p'}{q'} \in F_{n-1}$ we have $
\frac{p'}{q'}\oplus \frac{p^1}{q^1}= \frac{p^0}{q^0}$ and
$\frac{p'}{q'}, \frac{p^0}{q^0}$ lie on the same side from $x$. So
 (as it was mentioned above in the beginning of Section 6) from the Farey tree construction we see that $\xi^1 =
\frac{p_{t-1}}{q_{t-1}}$ is also a convergent fraction to $x$. So
\begin{equation}
\delta <\frac{1}{q_tq_{t-1}}. \label{case1}
\end{equation}
 We should consider two subcases: $a_{t+1} = 1$ ({\bf case 1.1}) and $a_{t+1} > 1$ ({\bf case 1.2}).

{\bf Case 1.1.} Here we have $a_{t+1} = 1$ . So $\xi = \frac{p_t+ p_{t-1}}{q_t+ q_{t-1}} =\frac{p_{t+1}}{q_{t+1}}$. Now we must look for the
natural number $z$ defined above. It may happen that $z=1$ ({\bf case 1.1.1}) or $z>1$ ({\bf case 1.1.2}).

{\bf Case 1.1.1.} Here from (\ref{w2w1}) we see that $ ?(x+\delta) -?(x) \ge \frac{1}{2^{n+2}}$ and together with (\ref{case1}) this leads to
the inequality
\begin{equation}
\frac{?(x+\delta) -?(x)}{\delta} \ge \frac{q_tq_{t-1}}{2^{n+2}} = \frac{q_tq_{t-1}}{2^{a_1+...+a_t+2}}
 . \label{case111}
\end{equation}

 {\bf Case 1.1.2.} Now $z >1$. From the Farey tree constructing process we see that $ z\le a_{t+2}$. As
$(x,x+\delta) \subset (\xi_{- -},\xi_{++})$ we see that
$$
\delta \le (\xi_{++}-\xi) +(\xi -  \xi_{- -})= \frac{1}{q_t+q_{t-1}}\left( \frac{1}{(z-1)(q_t+q_{t-1})+q_{t-1}} +
\frac{1}{(z-1)(q_t+q_{t-1})+q_{t}}\right)\le \frac{2}{(z-1)q_{t+1}^2}.
$$
Here we use the fact that  fractions $ \xi_{- -},\xi , \xi_{++}$ are successive elements from $F_{n+z}$ and $$|\xi_{++}-\xi| =
\frac{1}{(q_t+q_{t-1})((z-1)(q_t+q_{t-1})+q_{t-1})}, |\xi_{--}-\xi| = \frac{1}{(q_t+q_{t-1})((z-1)(q_t+q_{t-1})+q_{t})}. $$

 We apply (\ref{w2w1}) to get
$$
\frac{?(x+\delta) -?(x)}{\delta} \ge \frac{(z-1)q_{t+1}^2}{2^{a_1+...+a_t+z+2}} \ge \frac{zq_{t+1}^2}{2^{a_1+...+a_{t+1}+z+3}}.
$$
As $ 2\le z \le a_{t+2}$, $a_{t+2}q_{t+1} \ge  q_{t+2}/2$ and the function $ \frac{z}{2^z}$ decreases we see that in the case 1.1.2 the
following inequality is valid:
\begin{equation}
\frac{?(x+\delta) -?(x)}{\delta} \ge    \frac{q_{t+2} q_{t+1} }{2^{a_1+...+a_{t+1}+a_{t+2}+4}}
 . \label{case112}
\end{equation}

{\bf Case 1.2.} We have $ a_{t+1} \ge 2$ so $ z=1$ and $\xi_- = \frac{2p_t+p_{t-1}}{2q_t+q_{t-1}}\in (x,\xi )$. Now from (\ref{case1}) and
(\ref{w2w1}) we deduce that
\begin{equation}
\frac{?(x+\delta) -?(x)}{\delta} \ge    \frac{q_{t} q_{t-1} }{2^{n+2}} = \frac{q_{t} q_{t-1} }{2^{a_1+...+a_t+2}}
 . \label{case12}
\end{equation}

{\bf Case 2.} In this case we have $q_t = q^0<q^1$. We see that
$\xi^1   \in F_{ n}\setminus  F_{ n-1}$ and $\xi^0 =
\frac{p_t}{q_t} \in F_{n-1}$. So in the case 2 we see that
$n=a_1+...+a_t+l$ with some $l < a_{t+1}$.
 We can establish the inequality
\begin{equation}
\delta <\frac{1}{q_tq^1}. \label{case2}
\end{equation}
Consider the convergent  $\frac{p_{t-1}}{q_{t-1}}$. The fraction $\xi^1$ is an intermediate fraction to $x$. It lies between $ \xi$ and $\xi  '=
\frac{p_{t-1}}{q_{t-1}}$. Define $l$ from the condition $\xi^1 = \xi '\underbrace{\oplus \xi^0 \oplus ... \oplus \xi^0}_{l}$. So $ q^1=
q_{t-1}+lq_t$.   Moreover we see that
 $ a_{t+1} \ge l+1\ge 2$.
 It may happen that $\frac{p_{t+1}}{q_{t+1}} = \frac{p_t+p^1}{q_t+q^1} = \xi$
 ({\bf case 2.1}) or  $\frac{p_{t+1}}{q_{t+1}}\neq \xi$ ({\bf case 2.2}).

{\bf Case 2.1} Let $\frac{p_{t+1}}{q_{t+1}}  = \xi$. Then
$n=a_1+...+a_t+a_{t+1} -1$ and $ l = a_{t+1}-1$. It may happen
that $z=1$ ({\bf case 2.1.1}) or $z>1$ ({\bf case 2.1.2}).

{\bf Case 2.1.1.} Note that $ q^1+q_t = q_{t+1}=a_{t+1}q_t + q_{t-1}$ and so $ q^1 = q_{t+1} -q_t \ge (a_{t+1}-1) q_t+ q_{t-1} \ge q_{t+1}/2 $
as $ a_{t+1}\ge 2$.
 From (\ref{case2}) and (\ref{w2w1}) we see that
\begin{equation}
\frac{?(x+\delta) -?(x)}{\delta} \ge    \frac{q_{t} q^1 }{2^{n+2}} \ge \frac{q_{t} q_{t+1} }{2^{a_1+...+a_t+l+3}}\ge \frac{q_{t} q_{t+1}
}{2^{a_1+...+a_t+ a_{t+1}+2}}
 . \label{case211}
\end{equation}

{\bf Case 2.1.2.} Here we have $ 2\le z \le a_{t+2}+1$. As in the case 1.1.2 we have
$$
\delta \le   (\xi_{++}-\xi) +(\xi -  \xi_{- -})=   \frac{1}{q_{t+1}} \left( \frac{1}{(z-1)q_{t+1}+q_t}+ \frac{1}{(z-1)q_{t+1}+q^1}\right) \le
\frac{2}{(z-1)q_{t+1}^2}\le \frac{4}{zq_{t+1}^2} .
$$
From this inequality and (\ref{case2}) we see that
\begin{equation} \frac{?(x+\delta) -?(x)}{\delta} \ge    \frac{zq_{t+1}^2 }{2^{n+z+3}}\ge
\frac{(a_{t+2}+1)q_{t+1}^2 }{2^{n+a_{t+2}+3}} \ge
\frac{(a_{t+2}+1)q_{t+1}^2 }{2^{a_1 + ... + a_{t+1}+a_{t+2}+2}}
\ge
\frac{q_{t+2}q_{t+1}}{2^{a_1+...+a_{t+1}+a_{t+2}+2}}.
.
\label{case212}
\end{equation}

{\bf Case 2.2.} We have $\frac{p_{t+1}}{q_{t+1}}\neq \frac{p_t+p^1}{q_t+q^1}$. So in this case $ \frac{2p_t+p^1}{2q_t+q^1} =\xi_-\in (x, \xi )$
and $z =1$. We see from (\ref{case2}) and (\ref{w2w1}) that
$$
 \frac{?(x+\delta) -?(x)}{\delta} \ge \frac{q_tq^1}{2^{n+2}}.
$$
As $\frac{p^1}{q^1} = \frac{lp_t+p_{t-1}}{lq_t+q_{t-1}}$  and $ n = a_1+..+a_t + l$ we see that
\begin{equation}
 \frac{?(x+\delta) -?(x)}{\delta} \ge \frac{lq_t^2}{2^{a_1+...+a_t + l+2}}\ge
\frac{q_{t+1}q_t}{2^{a_1+...+a_t +a_{t+1}+ 3}}
 .
\label{case22}
\end{equation}

Now we pick together the results of the cases 1.1.1, 1.1.2, 1.2, 2.1.1, 2.1.2, 2.2 (that is the inequalities
(\ref{case111},\ref{case112},\ref{case12},\ref{case211},\ref{case212},\ref{case22})) to get the statement of the lemma. Lemma 1 is proved.

 {\bf Lemma 2.}\,\,\,{\it For
an irrational $x\in (0,1)$ and  for  $\delta  $ small enough  there exists a natural
  $ r = r (x,\delta )$ and  such that
\begin{equation}
\frac{?(x+\delta )-?(x)}{\delta} \le
   \frac{ q_{r+1}^2}{2^{a_1+...+a_{r+1}  -2}} \label{DDDD}
\end{equation}
}

{\bf Proof.}\,\, It is enough to prove Lemma 1 for positive  $\delta$.  We use the notation $\xi^0,=\frac{p_t}{q_{t}}, \xi^1
=\frac{p^1}{q^1},\xi, \xi_-,\xi_+,\xi_{--},\xi_{++}$ and $z$ from the proof of Lemma 1.
 As points $\xi_{--}<\xi_-<\xi <\xi_+<\xi_{++}$
are successive points from $F_{n+z+1}$
 and $?(x)$ increases, we have
\begin{equation}
  ?(x+\delta ) - ?(x) \le ?(\xi_{++}) - ?(\xi_{--}) = \frac{4}{2^{n+z+1}}.
  \label{w2}
\end{equation}

Consider two cases:   $\xi_- \in (x,\xi )$ ({\bf case 1}) and  $\xi_- \not\in (x,\xi )$ but then $\xi_{+} \in (\xi , x+\delta )$ ({\bf csae 2}).

{\bf Case 1.}\,\, Here we have $\delta > \xi - \xi_-$. It may happen that $z+1$  ({\bf case 1.1}) or $z>1$ ({\bf case 1.2}).

{\bf Case 1.1.}\,\, Let
 $ z=1$. Then $\xi_- =p/q,
q = z_* q_t + q_{t-1} \le q_{t+1}, 1\le z_* \le a_{t+1},
  \xi =(p - p_t)/(q-q_t)$,
$n+2 = a_1+...+a_{t}+z_*   $
  and
$$
\delta > \xi -  \xi_-\ge \frac{1}{(q-q_t)q} \ge \frac{1}{(z_*+1)^2q_{t}^2}.
$$
We take into account (\ref{w2}) to see that
\begin{equation}
\frac{?(x+\delta ) - ?(x) }{\delta} \le \frac{(z_*+1)^2q_t^2}{2^{n+2}}= \frac{(z_*+1)^2q_t^2}{2^{a_1+...+a_{t}+z_*}}\le
\frac{q_t^2}{2^{a_1+...+a_{t}-2}}.
 \label{26}
\end{equation}

{\bf Case 1.2.}\,\, If $ z > 1$
 then $\xi =p_{t+1}/q_{t+1}$, $ \xi_{--} = p_{t+2}/q_{t+2}, z =
a_{t+2}+1$, $n+1 = a_1+...+a_{t+1}$
 and
$$
\delta >  \xi -  \xi_-\ge \frac{1}{(zq_{t+1}+q_t)q_{t+1}} \ge \frac{1}{(z+1)q_{t+1}^2}.
 $$
From (\ref{w2}) we see that
\begin{equation}
\frac{?(x+\delta ) - ?(x) }{\delta} \le \frac{(z+1)q_{t+1}^2}{2^{n+z-1}}= \frac{(z+1)q_{t+1}^2}{2^{a_1+...+a_{t+1}+z-2}}\le
\frac{q_{t+1}^2}{2^{a_1+...+a_{t+1}-2}}.
 \label{27}
\end{equation}

{\bf Case 2.}\,\, Here we have $z\le a_{t+2},\xi = p_{t+1}/q_{t+1}, n+1 = a_1+...+ a_{t+1}$. Now we deduce
$$
\delta > \xi_+-\xi \ge \frac{1}{(zq_{t+1}+q^1)q_{t+1}}\ge \frac{1}{(z+1)q_{t+1}^2}
$$
(remind that $q^1 < q_{t+1} $ is the denominator of $\xi^1$). From (\ref{w2}) we see  that
\begin{equation}
\frac{?(x+\delta ) - ?(x) }{\delta} \le \frac{(z+1)q_{t+1}^2}{2^{n+z-1}}= \frac{(z+1)q_{t+1}^2}{2^{a_1+...+a_{t+1}+z-2}}\le
\frac{q_{t+1}^2}{2^{a_1+...+a_{t+1}-2}}.
 \label{28}
\end{equation}
In any case (\ref{DDDD}) follows from (\ref{26}),(\ref{27}) or (\ref{28}). Lemma 2 is proved.

{\bf 7.\,\,\, More notation.} In the sequel small letters $a,b,c$ will be used for natural numbers. By capital letters $A,B,C,... $ we denote
finite sequences of natural numbers. For
$$ A = a_1,a_2,...
, a_{t -1}, a_t
$$
we use the notation
$$\overrightarrow{A} = a_1,a_2,...,
a_{t -1}, a_t,\,\,\, \overleftarrow{A}= a_t,a_{t -1}, ..., a_2,a_1 . $$ If a sequence $A=a_1,...,a_t $ appears in a continuant expression or in
a continued fraction expression it means that we should replace it by the consecutive block $a_1,...,a_t $ of natural numbers $a_j$. For example
for $A=a_1,...,a_t $ ; $B =b_1,...,b_i$; $X = x_1,...,x_j $ and $Y = y_1,...,y_k $ we have
$$ \langle a,b,A,B,c,\overleftarrow{X},\overrightarrow{Y}\rangle = \langle
a,b,a_1,...,a_t,b_1,...,b_i,c, x_j,...,x_1,y_1,...,y_k\rangle
$$
and
$$ [  a,b,A,B,c,\overleftarrow{X},\overrightarrow{Y}]=[
a,b,a_1,...,a_t,b_1,...,b_i,c, x_j,...,x_1,y_1,...,y_k].
$$
For a sequence $A = a_1,
 a_2,... , a_{t -1}, a_t$ we define
$$
A_- = a_2,... , a_{t -1}, a_t,\,\,\, A^- = a_1,a_2,... , a_{t -1}.
$$
So
$$
[\overrightarrow{A}] =\frac{\langle A_- \rangle}{\langle A  \rangle},\,\,\,\,\, [\overleftarrow{{A}}] =\frac{\langle A^- \rangle}{\langle A
\rangle}
$$
and if $ x=[a_1,a_2,...,a_t,...]$ then
$$
\frac{p_t}{q_t} = \frac{p_t(x)}{q_t(x)}= \frac{\langle a_2,...,a_t \rangle}{\langle a_1,...,a_t  \rangle},\,\,\,\,\,\, \frac{q_{t-1}}{q_t} =
\frac{q_{t-1}(x)}{q_t(x)}= \frac{\langle a_1,...,a_{t-1} \rangle}{\langle a_1,...,a_t  \rangle} .
$$
Also in the case  when $A$ is an empty sequence we put $\langle A\rangle = 1$ and $\langle A_-\rangle = \langle A^-\rangle =[A] = 0$.

 This notation is convenient to work with some identities involving continuants and continued fractions. We
make use of the well known (see \cite{knuth})
 identity
 \begin{equation}
\langle X,Y\rangle = \langle X\rangle \cdot \langle Y\rangle +
 \langle X^-\rangle \cdot \langle Y_-\rangle
= \langle X\rangle \cdot \langle Y\rangle\cdot ( 1 +[\overleftarrow{X}]\cdot[\overrightarrow{Y}]).
 \label{knut}
 \end{equation}

 {\bf 8.\,\,\, Inequalities with  continuants: unit variation.}

 Let sequences $A,B,C$ be fixed.
Consider a natural number $\tau \ge 2$. For $a \in \mathbb{N},
a\le \tau -1$ we put $b = \tau - a$. So
$$a+b=\tau.$$
We define function
$$
F(a) =   F_{A,B,C;\tau }(a) = \langle A,a,C,\\b,B\rangle = \langle A,a,C,\tau - a,B\rangle .
$$
{\bf Lemma 3.}\,\,\,{\it Let $a,b\neq 1 $. Let $ n\in \mathbb{N} $ be an upper bound for the
maximal element among all the sequences $A,B,C$ and integers
$a,b$. Suppose $ [\overleftarrow{A}]\neq 0,1$.

 1.\, Let $a\neq b$. Then
$$
F(a) \ge \left( 1+\frac{1}{16(n+2)^3}\right)   F(a+ {\rm sign}((a-b)).
$$

2.\, Let $a=b$. Then
$$
F(a) \ge \left( 1+\frac{1}{16(n+2)^2}\right)  \min \{ F(a-1), F(a+1)\}.
$$

3.\,Let $a=b$,\,\, $ [\overrightarrow{C}]=[\overleftarrow{C}]$. Then
$$
F(a) \ge \left( 1+\frac{1}{16(n+2)^3}\right)  \max \{ F(a-1), F(a+1)\}.
$$
  }

{\bf Corollary 2.}\,\,\,{\it Let $a,b\neq 1 $. Let $ n $ be an upper bound for the maximal element among all the sequences $A,B,C$   and numbers $a,b$. Suppose $
[\overleftarrow{A}]\neq 0,1$. Then
$$
F(a) \ge \left( 1+\frac{1}{16(n+2)^3}\right)  \min \{ F(a-1), F(a+1)\}.
$$}

{\bf Definition 1}.\,\, Let $a_i\pm 1\ge 1, a_j\mp 1\ge 1$. We define the procedure
$$
  a _1,...,a_i,...,a_j,...,a_t  \,\,\,\,\mapsto\,\,\,   a_1,...,a_i\pm 1,...,a_j\mp 1,...,a_t  ,\,\,\, i<j
$$
as a {\it unit variation} of the  sequence $  a_1,...,a_i,...,a_j,...,a_t $.

{\bf Remark.} \,\, Let a sequence $a_1^*,..., a_t^*$ is obtained by a unit variation from the sequence
$a_1,...,a_t$. Then for any $\gamma \le t$ one has
$$\frac{1}{2} \langle a_1,...,a_\gamma\rangle \le
\langle a_1^*,...,a_\gamma^*\rangle \le 2\langle a_1,...,a_\gamma\rangle .
$$
(This Remark follows from formula (\ref{knut}).)

Proof of Lemma 3. We apply (\ref{knut}) to see that $F(a)$ is a
quadratic polynomial in $a$.
\begin{equation}
F(a) = \langle A,a \rangle \langle C \rangle \langle \tau - a, B\rangle + \langle A \rangle \langle C_- \rangle \langle \tau - a, B\rangle+
\langle A,a \rangle \langle C^- \rangle \langle B\rangle + W, \label{ww}
\end{equation}
 where $W$ does not depend on $a$.
  Note that
$$
\langle A,a \rangle \langle C \rangle \langle \tau - a, B\rangle = \langle A\rangle \langle B \rangle \langle C \rangle (a+
[\overleftarrow{A}])( \tau - a +[\overrightarrow{B}]),
$$
$$
\langle A\rangle \langle C_- \rangle \langle \tau - a, B\rangle = \langle A\rangle \langle B \rangle \langle C \rangle [ \overrightarrow{C}](
\tau - a +[ \overrightarrow{B}]),
$$
$$
\langle A, a \rangle \langle C^- \rangle \langle B\rangle= \langle A\rangle \langle B \rangle \langle C \rangle  [\overleftarrow{C}]( \ a +[
\overleftarrow{A}]).
$$
Put
$$
\theta = [\overrightarrow{B}]-[\overleftarrow{A}]-[\overrightarrow{C}]+[\overleftarrow{C}].
$$
 Now (\ref{ww}) changes to
\begin{equation}
F(a) =  \langle A\rangle \langle B \rangle \langle C \rangle((\tau +\theta)a-a^2) + W_1, \label{www}
\end{equation}
 where $W_1$ does not depend on $a$.

We want to obtain an upper bound for $|\theta |$. We see that
$[\overrightarrow{B}],[\overleftarrow{A}],[\overrightarrow{C}],[\overleftarrow{C}]\in [0,1]$. Moreover   $n$ is an upper bound for the elements
of the sequence $A$. So from the condition $ [\overleftarrow{A}]\neq 0,1$ we have $ \frac{1}{n+1}\le [\overleftarrow{A}]\le 1-\frac{1}{n+2}$.

 Note that in the case $
[\overrightarrow{C}]=[\overleftarrow{C}]$ we have
\begin{equation}|\theta |\le 1- \frac{1}{n+2}.
\label{320}
\end{equation}
 In the case   $ [\overrightarrow{C}]\neq [\overleftarrow{C}]$ we have the inequality
\begin{equation}
|\theta | \le 2 - \frac{1}{n+2}. \label{32}
\end{equation}
From (\ref{knut}) we see that
\begin{equation}
\langle A, a\pm 1, C,b\mp 1,B\rangle\le 16 \langle A\rangle \langle B \rangle \langle C \rangle (a\pm 1)(b\mp 1)\le 16(n+1)^2\langle A\rangle
\langle B \rangle \langle C \rangle \label{ind1}
\end{equation}
 as $ 2\le a,b\le n$. Now
\begin{equation}
 \langle A\rangle \langle B \rangle \langle C \rangle \ge \frac{\langle A, a\pm 1, C,b\mp 1,B\rangle}{16(n+1)^2}=
\frac{F(a\pm 1)}{16(n+1)^2}
 .
\label{wwwww}
\end{equation}

From (\ref{www}) with $\varepsilon \in \{-1,+1\}$ we have

\begin{equation}
F(a) = F(a +\varepsilon)+ \langle A\rangle \langle B \rangle \langle C \rangle(\varepsilon(a-b-\theta) +1)
 .\label{ind2}
\end{equation}
 Now we substitute (\ref{wwwww}) into the last inequality and obtain
$$
F(a) \ge \left( 1+\frac{\varepsilon(a-b-\theta) +1}{16(n+2)^2}\right)  F(a+\varepsilon ) .
$$
Let $a\neq b$ and $\varepsilon = {\rm sign}(a-b)$. Then by (\ref{32}) we have
$$
\varepsilon(a-b-\theta) +1 = |a-b|+1 \pm \theta \ge 2 - |\theta |\ge  \frac{1}{n+2}.
$$
So we have
\begin{equation}
F(a) \ge \left( 1+\frac{ 1}{16(n+2)^3}\right)  F(a+{\rm sign}(a-b) ) \label{ind3}
\end{equation}
and the first statement of Lemma 3 is proved.

 Let $a=b$. Then
\begin{equation}
F(a) \ge \left( 1+\frac{-\varepsilon\theta +1}{16(n+2)^2}\right)
F(a+\varepsilon ) \label{OO}
\end{equation}
 Now we choose $\varepsilon $ to
satisfy $ -\varepsilon \theta \ge 0$ and obtain the second
statement of Lemma 3.

Now we must consider the case when in addition to $a=b$ we suppose $ [\overrightarrow{C}]=[\overleftarrow{C}]$. In this case from (\ref{OO}) and
(\ref{320}) it follows that
 $$
F(a) \ge \left( 1+\frac{1-|\theta|}{16(n+2)^2}\right)
F(a+\varepsilon )\ge  \left( 1+\frac{1}{16(n+2)^3}\right)
F(a+\varepsilon).
$$
 Lemma 3 is proved.

We need an improvement of Lemma 3 in the case $ a\le 4 < 5\le b$.

 {\bf Lemma
4.}\,\,\,{\it Let $ a\le 4 < 5\le b= \tau - a,\,\, a+2\le b$. Then
$$
 F(a+1)
 \ge \frac{385}{384} F(a).
$$
} Proof. We follow the arguments of the proof of Lemma 3. As $ a \le 4$ instead of the inequality (\ref{ind1}) we see that
$$
F(a) \le 16  \langle A\rangle \langle B \rangle \langle C \rangle ab \le 64 b \langle A\rangle \langle B \rangle \langle C \rangle.
$$
So instead of (\ref{ind2},\ref{ind3})($\varepsilon = -1,\theta \ge -2$) we get
$$
F(a+1) = F(a)+  \langle A\rangle \langle B \rangle \langle C \rangle (b-a+\theta+1)\ge F(a)\left( 1+\frac{b-\min (5,b-1)}{64b}\right) \ge \frac{385}{384} F(a)
$$ as $b \ge 5$. Lemma is proved.

 Define
$$
U(t,s) = \{ (a_1,a_2,...,a_t) :\,\, a_j \in \mathbb{N},
a_1+a_2+...+a_t=s\},
$$
$$
U_n(t,s) = \{ (a_1,a_2,...,a_t) :\,\, a_j \in \mathbb{N}, \, a_j \le n,\,\,  a_1+a_2+...+a_t=s\}.
$$

 {\bf Lemma
5.}\,\,\,{\it

1. The minimal value of the continuant
$$\langle 1,1,a_1,a_2,...,a_t\rangle
=
\langle 2,a_1,a_2,...,a_t\rangle,\,\,\,(a_1,a_2,...,a_t)\in
U_n(t,s)
$$
attains at the sequence $Q =a_1,a_2,...,a_t$ such that not more than one element $a_j$ differs from $1$ and from $n$.

2. The maximal value of the continuant
$$
\langle 1,1,a_1,a_2,...,a_t\rangle = \langle
2,a_1,a_2,...,a_t\rangle,\,\,\, (a_1,a_2.,..,a_t)\in U(t,s)
$$
attains at the sequence $P=a_1,a_2,...,a_t$ such that for every
$i,j \in \{1,2,...,t\}$ one has $|a_i-a_j|\le 1$.}

Proof.

1. Suppose the minimum attains at a sequence with some elements $
a_i, a_j \not \in \{1,n\}, i< j$. We put $ a= a_i, b=a_j, \tau =a+b$ and apply
Corollary 2 to see that
$$
\langle 2,a_1,...,a_i,...,a_j,...,a_t\rangle > \min \{ \langle 2,a_1,...,a_i-1,...,a_j+1,...,a_t\rangle, \langle
2,a_1,...,a_i+1,...,a_j-1,...,a_t\rangle\}
$$
and this is a contradiction as
$$
(a_1,...,a_i-1,...,a_j+1,...,a_t),(
a_1,...,a_i+1,...,a_j-1,...,a_t) \in U_n(t,s).
$$

2. Assume that the conclusion is not true. Then we take one of the shortest subsequence  $a_i,a_{i+1},...,a_j$ of $P$ such that $|a_i-a_j| \ge
2.$ Suppose $ a_i >a_j$ without loss of generality. Put $a = a_i -1, b=a_j+1, \tau = a+b$.

If $a> b$ we have $a,b\ge 2, \,{\rm sign}(a-b) =+1.$ So we apply
statement 1 of Lemma 3 and obtain
$$
\langle 2, a_1,...,a_i-1,...,a_j+1,...,a_t\rangle
> \langle 2, a_1,...,a_i,...,a_j,...,a_t\rangle
.
$$ So $
\langle 2, a_1,...,a_i,...,a_j,...,a_t\rangle$ is not maximal.

If $a=b$ we
 make use of the fact that the considered subsequence $a_i,a_{i+1},...,a_j$ is the shortest one. We have
$$
a_{i+1}=a_{i+2} =...=a_{j-1} = a_i-1=a_j+1.\
$$
So for the sequence $C = a_{i+1},a_{i+2},...,a_{j-1}$ we have $
[\overrightarrow{C}]=[\overleftarrow{C}]$. Moreover $A=2,a_1,...,
a_{i -1}$ is such that $[\overleftarrow{A}] \neq 0,1 $. Now we can
apply statement 3 of Lemma 3 and it follows that
$$
\langle 2,a_1,...,a_i-1,...,a_j+1,...,a_t\rangle
>
  \langle 2,a_1,...,a_i,...,a_j,...,a_t\rangle .
 $$
 We see again that    $
\langle 2,a_1,...,a_i,...,a_j,...,a_t\rangle$ is not maximal.

Lemma 5 is proved.

We need the following supplement to Lemma 5.

 {\bf Lemma
6.}\,\,\,{\it

1. For any sequence $(a_1,a_2,...,a_t)\in U_n(t,s)$ one has
\begin{equation}\langle 1,1,a_1,a_2,...,a_t\rangle
 \ge \left( 1+\frac{1}{16(n+2)^3}\right)^{\frac{2\sigma _1  -n}{4}} \times
\min_{(b_1,b_2,...,b_t)\in U_n(t,s)} \langle
1,1,b_1,b_2,...,b_t\rangle, \label{ggg}
\end{equation}
 where
\begin{equation}
 \sigma_1=
 \sigma_1 (a_1,a_2,...,a_t)=
 \sum_{j=1}^t \min (a_j-1, n
-a_j) .
\label{1w1}
\end{equation}

2. Let the maximal continuant  with elements from the set  $U(t,s)$ contains elements $a =4 ,b= 5$. Then for any sequence
$(a_1,a_2,...,a_t)\in U(t,s)$ one has
\begin{equation}\langle 1,1,a_1,a_2,...,a_t\rangle
 \le \left( \frac{385}{384}\right)^{-\frac{\sigma_2}{2}} \times
\max_{(b_1,b_2,...,b_t)\in U(t,s)} \langle
1,1,b_1,b_2,...,b_t\rangle,
 \label{gggg}
\end{equation}
where
 \begin{equation}
 \sigma_2= \sigma_{2,t}=
 \sigma_2 (a_1,a_2,...,a_t)=
\sum_{j=1}^t \min(|a_j-4|,|a_j-5|)
 \label{mm1}
 \end{equation}
}

Proof.

To prove Lemma 6 we observe that we can obtain the sequence with relatively  minimal (maximal) continuant from any given sequence $(a_1,a_2,...,a_t)\in
U_n(t,s)$
(or $\in U(t,s)$)
by successive applications of the unit variation procedure (see Definition 1). Here we must use only those unit variations for
 which the value of the continuant
decreases (increases).
 From
the proof of  Lemma 1 we see that it is possible to do indeed.

Now we prove the statement 1. Note that for $1<z<n$ the following inequality is valid:
$$
| \min(z-1, n-z)- \min(z-2, n-z+1)|\le 1.
$$
When a pair of elements $a_i,a_j$ is replaced  by the pair $a_i-1,a_j+1$ we have
$$
|\sigma_1 (a_1,..., a_i-1,...,a_j+1, ...,a_t)-\sigma_1 (a_1,..., a_i,...,a_j, ...,a_t)|\le 2.
$$
After we have made all necessary unit variations we come to the minimal continuant $\langle 1,1,a_1^*,a_2^*,...,a_t^*\rangle$ where all elements
but one are from the set $\{1,n\} $. This exceptional  element (if exists) must be greater than $1$ and less than $n$.
 Hence
$$
0\le  \sigma_1 (a_1^*,a_2^*,...,a_t^*) \le n/2.
$$
Now we see that the number $U$ of unit variations (used to get from the initial  continuant the minimal one) is not less than the number $U'$ of
unit variations for  which the sum $ \sigma_1 (a_1,... ,a_t)$ strictly decreases. So
$$
U\ge U'\ge
 \frac{\sigma _1 (a_1,...,a_t) -n/2}{2}
 .$$
 Each unit variation enlarges the continuant by the factor $\left( 1+\frac{1}{16(n+2)^3}\right) $ as it was shown in
  Lemma 3.
 The statement (\ref{ggg}) follows.

Now we prove the statement 2. Note that
$$
| \min(|z-4|, |z-5|)- \min(|z-4-1|, |z-5-1|)|\le 1.
$$
 When a pair of elements $a_i,a_j$   is replaced  by the pair $a_i+1,a_j-1$ we have
$$
|\sigma_2 (a_1,..., a_i+1,...,a_j-1, ...,a_t)- \sigma_2 (a_1,..., a_i,...,a_j, ...,a_t)|\le 2.
$$
For the maximal continuant $\langle 1,1,a_1^{**},a_2^{**},...,a_t^{**}\rangle$ we have obviously
$$
\sigma_2 (a_1^{**},a_2^{**},...,a_t^{**})=0.$$ So the number of the unit variations used in the process of getting the maximal continuant is not
less than $ \sigma_2 (a_1, ...,a_t)/2$.

Note that we know that the maximal continuant contains of digits $4,5$ only. So given sequence $a_1,...,a_t$ from $U(t,s)$ we may assume that
some of digits $a_j$ are $\le 4$ meanwhile some of digits $a_j$ are $ \ge 5$. Moreover to obtain  the sequence with the maximal value of the
continuant from the sequence $a_1,...,a_t$ we may use the unit variation procedures with replacing digits $(a,b)\mapsto (a+1, b-1),\,\,\,
a+1\le b-1,\,\,\,
a\le
4<5\le b$ only. Now (\ref{gggg}) follows from Lemma 4.

Lemma 6 is proved.

 {\bf 9.\,\,\, Inequalities with continuants: substitutions.}
Put $a_0 = 2$.
 In this section we consider three sequences $A=a_0,a_1,...,a_w;\,\,\,B=b_1,...,b_r$
and $C=c_1,c_2,...,c_k $. Define
\begin{equation}
K =\langle A,C, B \rangle, \,\,\, \alpha = \alpha (A,B,C)=
(c_1-c_k)([\overleftarrow{A}]-[\overrightarrow{B}]), \,\,\, \Psi =
\Psi(A,B,C)= \langle A,\overleftarrow{C} , B \rangle .
\label{alef}
\end{equation}

 We
suppose $ n$ to be the maximal element among the elements of the
sequences $A,B ,C$.

 {\bf Lemma 7.}\,\,\,{\it Let $\alpha  = \alpha (A,B,C) \neq 0$. Let $w,r\ge 1 $ and $ a_w\neq b_1$ or $r\ge 2,\,\, a_{w
-1} \neq b_{2} $. Then
$$
1+  \frac{1}{4(n+2)^8} \le \left(\frac{\Psi }{K}\right)^{{\rm
sign}(\alpha )} \le 2.
$$}

{\bf Remark.}\,\, Obviously Lemma 7 remains  true in the case when $B$ is an empty sequence and $A,C$ are nonempty as formally $[B] =0$.

 Proof.\,\, It follows from (\ref{knut}) and
$\langle \overleftarrow{X}\rangle=
\langle\overrightarrow{X}\rangle
 $ that
$$
K = \langle A\rangle \langle C\rangle\langle B\rangle + \langle A^-\rangle \langle C_-\rangle\langle B\rangle + \langle A\rangle \langle
C^-\rangle\langle B_-\rangle+ \langle A^-\rangle \langle C^-_-\rangle\langle B_-\rangle
$$
and
$$
\Psi (A,B,C) = \langle A\rangle \langle C\rangle\langle B\rangle + \langle A^-\rangle \langle C^-\rangle\langle B\rangle + \langle A\rangle \langle
C_-\rangle\langle B_-\rangle+ \langle A^-\rangle \langle C^-_-\rangle\langle B_-\rangle
$$
that
\begin{equation}
K-\Psi  = \langle A\rangle \langle B\rangle\langle C\rangle (
[\overrightarrow{C}]-
[\overleftarrow{C}])([\overleftarrow{A}]-[\overrightarrow{B}]).
\label{pp}
\end{equation}
    As $ \alpha(A,B,C)\neq
0$ we see that ${\rm sign}(c_1-c_k) = -{\rm sign} ( [\overrightarrow{C}]- [\overleftarrow{C}])$. So
$$
{\rm sign} (\Psi  - K) = +{\rm sign} (\alpha(A,B,C)).
$$
We should take into account  inequalities
\begin{equation}
\frac{K}{4} \le \langle A\rangle \langle B\rangle\langle C\rangle \le K \label{ppp}
\end{equation}
\begin{equation}
\frac{\Psi }{4} \le \langle A\rangle \langle B\rangle\langle
C\rangle \le \Psi  \label{ppp1}
\end{equation}
 and
\begin{equation}
\frac{1}{(n+2)^3} \le |[\overrightarrow{C}]- [\overleftarrow{C}]|  \le 1,\,\,\, \frac{1}{(n+2)^5} \le
  |[\overleftarrow{A}]-[\overrightarrow{B}]|\le 1.
 \label{pppp}
\end{equation}
The upper bounds here are obvious. We give  our comments to the lower bound from (\ref{pppp}). As $c_1\neq c_k$ we have
$$
 |[\overrightarrow{C}]- [\overleftarrow{C}]|\ge
 \min_{1\le c\le n}
 \left(
\frac{1}{c+1}- \frac{1}{c+1+1/(n+1)}
 \right)=\frac{1}{(n+1)^2(n+1+1/(n+1))}\ge \frac{1}{(n+2)^3}.
 $$
If $ a_w\neq b_1$ by the similar reasons we have
$$
|[\overleftarrow{A}]-[\overrightarrow{B}]|\ge \frac{1}{(n+2)^3}.
$$
If $ a_w = b_1$ but $a_{w - 1} \neq b_{2}$ we see that
$$
|[\overleftarrow{A}]-[\overrightarrow{B}]|\ge \min_{1\le a_w, a_{w - 1} \le n }\left( \frac{1}{a_w+1/(a_{w - 1}+1/(n+1))}- \frac{1}{a_w+1/a_{w -
1}} \right)
>\frac{1}{(n+2)^5}
$$
and (\ref{pppp}) follows.

If ${\rm sign} (\alpha )=+1$ we deduce from (\ref{pp}) the
following inequality:
$$
1< \frac{\Psi }{K} = 1-\frac{\langle A\rangle \langle
B\rangle\langle C\rangle ( [\overrightarrow{C}]-
[\overleftarrow{C}])([\overleftarrow{A}]-[\overrightarrow{B}])}{K}.
$$
Now   from    the last formula and (\ref{ppp},\ref{pppp}) we have
$$
1
 + \frac{1}{4(n+2)^8}
< \frac{\Psi }{K} \le 2
$$
and Lemma 7 is proved in this case.

If ${\rm sign} (\alpha )=-1$   from (\ref{pp}) we have the
following inequality:
$$
1< \frac{K}{\Psi  } = 1+\frac{\langle A\rangle \langle
B\rangle\langle C\rangle ( [\overrightarrow{C}]-
[\overleftarrow{C}])([\overleftarrow{A}]-[\overrightarrow{B}])}{\Psi
(A,B,C)}.
$$
Now   from    the last formula and (\ref{ppp1},\ref{pppp}) we have
$$
1 + \frac{1}{4(n+2)^8}
 < \frac{K}{\Psi } \le 2
$$
and Lemma 7 is proved.

{\bf Remark.}\,\, Lemma 7 is close to the inequalities considered in \cite{knuth},\cite{knuth1}. Equalities similar to (\ref{pp}) were
considered in \cite{MOZ}.

We consider a sequence $D=d_1,...,d_t$ of the length $t$. Suppose
that $h_1<h_2<...<h_f=n= \max_{1\le i\le t} d_i$ are all  {\it different} elements of the
sequence $D$. Define
$$
r_\nu = r_\nu (D) = \# \{ i:\,\, d_i =\nu\}\ge 0 ,\,\,\, \nu =1,2,...,n.
$$
Of course $ r_\nu \ge 1$ iff $\nu = h_j$ for some $j$.
Let
$$
\pi = \left(
\begin{array}{cccc}
1& 2 &\cdots & t\cr \pi (1) &\pi (2)&\cdots & \pi(t)
\end{array}
\right)
$$
be a substitution.

The following lemma was announced in \cite{KAN} without a proof.
Here we give a complete proof.  More general setting was considered in \cite{knuth}.

{\bf Lemma 8.}\,\,\,{\it For a given $D$ the following equality is valid for the maximum over all substitutions:
$$
\max_\pi \langle 1,1,d_{\pi (1)},  d_{\pi (2)} ,\dots , d_{\pi (t)}
\rangle = \langle 1,1, \underbrace{1,...,1}_{r_1},
\underbrace{2,...,2}_{r_2},\dots ,
\underbrace{n,...,n}_{r_n}
 \rangle
.$$}

Proof. Consider any sequence $a_1,...,a_t$ of natural numbers such
that $ a_i >a_j$ with some $i<j$. Let  $ a_{i^*}, a_{j^*}$ be  two
elements of this sequence such that  $ a_{i^*} >a_{j^*}$ with
maximal value of the difference $ j^* - i^*$. Then
\begin{equation}
a_{i^*-1} \le a_{j^*} <a_{i^*} \le a_{j^*+1}.
 \label{tuda}
\end{equation}
(Of course it may happen that $i^*=1$ or $j^*=t$. Then  $ a_{i^*-1}= 1$ or the sequence $B$ below is an  empty sequence.)
  Put
$$
K = \langle 2, a_1,...,a_t\rangle = \langle A,C,B\rangle,\,\,\,
$$
where
$$
A = 1,1,a_1,...,a_{i^*-1};\,\,\,\, B= a_{j^*+1},...,a_t;\,\,\, C= a_{i^*},...,a_{j^*} .
$$
From (\ref{tuda}) it follows that for $\alpha (A,B,C)$  defined in
(\ref{alef}) we have
$$
\alpha (A,B,C) = (a_{i^*} -a_{j^*})
([\overleftarrow{A}]-[\overrightarrow{B}]) >0.
$$
Now from Lemma 7 we see that
$$
\langle A,\overleftarrow{C},B \rangle = \Psi (A,B,C) > K = \langle
A,C,B\rangle.
$$
But the sequence $ A,\overleftarrow{C},B$ may be obtained from the sequence $a_1,...,a_t$ by a certain substitution. So for any sequence
$a_1,...,a_t$ with an inversion  $ a_i >a_j, i<j$ we can find a permutation $\pi$ which enlarges the continuant $\langle a_1,...,a_t \rangle$.
Lemma 8 is proved.

We need a supplement to Lemma 8 in the case when the sequence $D$ consists of elements $4,5$ only.

{\bf Lemma 9.}\,\,\,{\it Let $D$ consists of elements $4,5$ only.
 Then
$$
  \langle 1,1,d_{1},  d_{2} ,\dots , d_{t}
\rangle \le  \left( 1+  \frac{1}{4\cdot 7^8}  \right)^{-\rho/2}
\langle 1,1, \underbrace{4,...,4}_{r_4},
\underbrace{5,...,5}_{r_5}
 \rangle
$$
where
\begin{equation}
\rho =\rho_t= \rho (d_{1},  d_{2} ,\dots , d_{t})= \sum_{j=1}^{t-1} |d_j -
d_{j+1}|- 1.
\label{mm2}
\end{equation}

}

Proof. By means of the substitutions from the proof of Lemma 8 one can  construct from the initial continuant $
  \langle 1,1,d_{1},  d_{2} ,\dots  ,d_{t}
\rangle $ the maximal one. After a substitution from the process initiated by the proof of Lemma 7 the sum
\begin{equation}
 \sum_{j=1}^{t-1} |d_j -
d_{j+1}| \label{deltiii}
\end{equation}
 decreases by 2.
For the maximal continuant the sum (\ref{deltiii}) is  not greater than $1$. So the number of the substitutions applied is  not less than $
\rho/2$. Each substitution enlarges the continuant by the factor $\left( 1+ \frac{1}{4\cdot 7^8}  \right)$ as it was shown in Lemma 7. Lemma is
proved.

{\bf Lemma 10.}\,\,\,{\it Let for the sequence $D$ under consideration we have $h_1=1$ and $r_1\ge t/2$. Then the following equality is valid for
the minimum over all substitutions:
$$
\min_\pi \langle 1,1,d_{\pi (1)},  d_{\pi (2)} ,\dots ,  d_{\pi (t)}
\rangle =
$$
$$
= \langle \underbrace{1,...,1}_{2+r_1-r_2-...-r_n}, \underbrace{2,1,2,1,...,2,1}_{2r_2\,\text{elements}},
\underbrace{3,1,3,1,...,3,1}_{2r_3\,\text{elements}},\dots , \underbrace{n,1,n,1,...,n,1}_{2r_n\,\text{elements}}
 \rangle
.$$}

Proof. First of all we prove that if  a substitution $a_1,...,a_t$
of the sequence $D$ give the minimal value of the continuant
$\langle 1,1, a_1,...,a_t \rangle $
 then among every two consecutive elements
$a_i,a_{i+1}$ at least one element is equal to $1$.

Suppose it is not so and we have two consecutive elements
$a_{i^*},a_{i^*+1} >1$. As $r_1> t/2$ we observe that there exists
a pair of consecutive elements $a_l, a_{l+1} = 1,1$. So
$$
1,1, a_1,...,a_t = 1,1, a_1,...,a_{l-1},1,1, a_{l+2},..., a_{i^*},a_{i^*+1},a_{i^*+2},... ,a_t
$$
or
$$
1,1, a_1,...,a_t = 1,1, a_1,..., a_{i^*},a_{i^*+1},a_{i^*+2},... , a_{l-1},1,1, a_{l+2},...,
  ,a_t.
$$
We consider only  the first opportunity as the second one is
similar. Put
$$
A = 1,1, a_1,...,a_{l-1},1;\,\,\,\, B =a_{i^*+1},a_{i^*+2},...
,a_t;\,\,\,\, C = 1, a_{l+2},..., a_{i^*}.
$$
For $\alpha (A,B,C)$ defined in (\ref{alef}) we have
$$
\alpha (A,B,C) = (1-
a_{i^*})([\overleftarrow{A}]-[\overrightarrow{B}]) <0.
$$
Hence
$$
\langle 1,1, a_1,...,a_t \rangle = \langle A,C,B\rangle > \Psi
(A,B,C).$$ So the permutation $a_1,...,a_t $ does not give the
minimal value of the considered continuant. We came to the
conclusion that the minimal continuant must be of the form
\begin{equation}
\langle 1,1,\underbrace{1,...,1}_u,a_1,\underbrace{1,...,1}_{v_1},a_2,\underbrace{1,...,1}_{v_2},...,a_{w-2},\underbrace{1,...,1}_{v_{w-2}},
 a_{w-1},\underbrace{1,...,1}_{
 v_{
w -1}}
 ,a_w,\underbrace{1,...,1}_{
 v_{
w }}\rangle,\,\,\, a_j\ge 2, \, j = 1,2,...,w, \label{mini}
\end{equation}
 where $ w =r_2+...+r_f,  v_1+...+v_w+u = r_1$, and $v_j \ge 1, 1\le j \le w - 1$.

Now we  shall show that $v_j=1$ for all $j \in \{1,...,w\}$.

First of all we shall show that $v_j\le 1$ for all $j$. Indeed if $v_j\ge 2$ for some j we put
$$
A = 1,1,\underbrace{1,...,1}_u;\,\,\,\, B =1, a_{j+1} ,\underbrace{1,...,1}_{
 v_{
j +1}}, ..., a_w,\underbrace{1,...,1}_{
 v_{
w }}
 ;\,\,\,\,
  C= a_1,\underbrace{1,...,1}_{v_1},a_2,\underbrace{1,...,1}_{v_2},...,
 a_j,\underbrace{1,...,1}_{
 v_{
j }-1}
$$
  Then
$$
\alpha (A,B,C) = (a_1-1)([\overleftarrow{A}]-[\overrightarrow{B}]) <0
$$
and by Lemma 7 we have
$$
\langle A,C,B \rangle > \Psi (A,B,C).
$$
So $ \langle A,C,B \rangle   $ is not minimal.

We show that $v_w\neq 0$. This statement is obvious in the case $w=1$. Suppose $ w\ge 2$. Then if $v_w=0$ put
$$
A = 1,1,\underbrace{1,...,1}_u,a_1;\,\,\,\,
  C=  1,a_2,1,...,a_{w-2},1, a_{w-1},1,a_{w}
$$
and $B$ define to be the empty sequence.
  Now
$$
\alpha (A,B,C) = (1-a_w)[\overleftarrow{A}] <0
$$
and  again by Lemma 7 we se that
 $ \langle A,C,B \rangle   $ is not minimal.

We have proven that the minimal continuant is of the form
\begin{equation}
\langle 1,1,\underbrace{1,...,1}_{r_1-r_2-...-r_n},a_1,1,a_2,1,...,a_{w-2},1,
 a_{w-1},1,a_w,1\rangle,\,\,\,
a_j\ge 2, \, j = 1,2,...,w,\,\,\,   w =r_2+...+r_n. \label{minii}
\end{equation}
 To finish the proof of Lemma 10 we must show that in the minimal continuant  $
a_i\le a_j$ for $i\le j$. Suppose that for some $i\le j$ one has $ a_i> a_j$. Let  $ a_{i^*}, a_{j^*}$ be  two elements of this sequence such
that  $ a_{i^*} >a_{j^*}$ with maximal value of the difference $ j^* - i^*$. Then as in the proof of Lemma 8 we have
$$
a_{i^*-1} \le a_{j^*} <a_{i^*} \le a_{j^*+1}.
$$
 We take
 $$
A = 1,1,a_1,1,...,a_{i^*-1},1;\,\,\,\, B= 1,a_{j^*+1},1,...,a_t;\,\,\, C= a_{i^*},1,a_{i^*+2}1,...,a_{j^*} .
$$
Then
$$
\alpha (A,B,C) = (a_{i^*}-a_{j^*})([\overleftarrow{A}]-[\overrightarrow{B}]) <0
$$
as $[\overleftarrow{A}]<[\overrightarrow{B}] $. Again from Lemma 7 we see that  $ \langle A,C,B \rangle   $ is not minimal. So we have proven
that $ a_i\le a_j$ for all $i\le j$. Lemma 10 follows.

We also need a supplement to Lemma 10 in the case when the sequence $d_1,...,d_t$ consists of elements $1,n$ only.
First of all given a sequence of partial quotients $d_1 , d_2 , ..., d_t$
we apply  certain consecutive substitutions to obtain a sequence $d_1^* , d_2^* , ..., d_t^*$
such that for any $i$  either $d_i^*$ or $d_{i+1}^* $ is equal to $1$.
We may take a certain type of this procedure to ensure that the sequence $d_1^* , d_2^* , ..., d_t^*$ depends on the initial sequence
 $d_1 , d_2 , ..., d_t$ only.

{\bf Lemma 11.}\,\,\,{\it Suppose that the conditions of Lemma 10 are satisfied. Moreover suppose that the sequence $d_1,...,d_t$ consists of
elements $1,n$ only. Then
 we have
$$
\langle 1,1,d_{1},  d_{2} ,\dots , d_{t} \rangle \ge \left(1+  \frac{1}{4(n+2)^8} \right)^\omega\,\,
 \min_\pi \langle 1,1,d_{\pi (1)},  d_{\pi (2)} ,\dots  , d_{\pi (t)} \rangle  .
$$
where
\begin{equation}
\omega = \omega (d_1 , d_2 , ..., d_t
 )=
  \sum_{j\le t-1: \,\,d_j \ge 2}\delta(d_j,d_{j+1})
+
\frac{1}{2}\sum_{j=1}^{t-2}
(1-\delta (d_j^*, d_{j+2}^*)) -1
  ,\,\,\,\, \delta (a,b) =
\begin{cases} 1,\,\, a=b,\cr
0,\,\, a\neq b.
\end{cases}
\label{2w2}
\end{equation}

}

Proof. We see from the proof of Lemma 10 that we can obtain the minimal continuant from a given continuant in the following manner. First of all
we apply substitutions to  ensure that among two successive elements one element is equal to 1.  For each permutation of such a kind the sum
\begin{equation} \sum_{j\le t-1: \,\,d_j \ge 2}\delta(d_j,d_{j+1})
\label{delti}
\end{equation}
 decreases by $1$
(remind that all $d_j\ge 2$ are equal to $n$). After several
substitutions we come to a continuant of the form (\ref{mini}) and
the sum (\ref{delti}) will be equal to zero. So we have used
exactly $\sum_{j\le t-1: \,\,d_j \ge 2}\delta(d_j,d_{j+1})$
substitutions on this stage. The sequence $d_1,...,d_t$ is
transformed now into the sequence $d_1^*,...,d_t^*$.

Then we use substitutions to pass from a continuant of the form (\ref{mini}) to a continuant of the form (\ref{minii}). During this process the
sum (\ref{delti}) does not change and remains equal lo zero. As for the sum
\begin{equation}
 \sum_{j=1}^{t-2} (1-\delta (d_j^*,d_{j+2}^* ))\label{deltii}
\end{equation}
we have the following observation. After each   permutation this sum now
  decreases
 by $ 2$. For the minimal continuant the sum (\ref{deltii}) is less or equal to $ 1$ (there is not more than  one nonzero summand corresponding to the
 first element which is not equal to $1$).
So we have used $\omega = \omega (d_1 , d_2 , ..., d_t
 )
$
permutations to pass from the initial continuant to the minimal one. Each permutation adds at least a factor $\left(1+ \frac{1}{4(n+2)^8}\right)$ by
Lemma 7.

Lemma is proved.

 {\bf 10.\,\,\, Some estimates.}\,\,\,
 Consider for  irrational $x$  the continued fraction expansion  $x=[a_1,a_2,...,a_t,...]$. For any naturals $k
, t $
 we define
 \begin{equation}
 w_k(
t) = \# \{ j\le t:\,\,\, a_j\ge k\}, \,\,\,\,\,\,
 r_k(
t) = w_{k}(t) -w_{k+1}(t)=
 \# \{ j\le t:\,\,\, a_j= k\}.
\label{WE}
\end{equation}

{\bf Lemma 12.}\,\,\,{\it Suppose that for some $t$ for an irrational number $x=[a_1,a_2,...a_t,...]$ we have $2r_1(t)\ge t$. Then for
$ n= n (t) =
\max_{j\le t} a_j$
 one has
\begin{equation}
\langle 2, a_1, ..., a_t \rangle \ge
\frac{\lambda_1^{t}}{10}\prod_{k=2}^{n} \left(
\frac{\mu_k}{\lambda_1^2} \right)^{r_k(t)}.
 \label{L7}
\end{equation}
where $\lambda_1$ and $\mu_k$ are defined in (\ref{lambdaj}), (\ref{muj}) respectively.}

As $ \lambda_1^2 = \mu_1$ and $w_{n+1} (t) = 0$ we see that
$$
\prod_{k=2}^{n} \left( \frac{\mu_k}{\lambda_1^2} \right)^{r_k(t)} = \prod_{k=2}^{n} \left( \frac{\mu_k}{\mu_{k-1}} \right)^{w_k(t)}.
$$
Obviously $\mu_{k+1} >\mu_k$.
 So  we immediately obtain

{\bf Corollary 3.}\,\,\,{\it Under conditions of Lemma 12 one has for any natural $N\le n$ the following inequality
$$
\langle 2, a_1, ..., a_t \rangle \ge
\frac{\lambda_1^{t}}{10}\prod_{k=2}^{N} \left(
\frac{\mu_k}{\mu_{k-1}} \right)^{w_k(t)}.
$$}

Proof of Lemma 12.\,\,\,
  By Lemma 10 we see that
 $$
 \langle 2, a_1, ..., a_t \rangle \ge
\langle \underbrace{1,...,1}_{1+t-2(r_2+...+r_{n})}, \underbrace{1,2,...,1,2}_{2r_2\,\text{elements}},...,
\underbrace{1,n,...,1,n}_{2r_{n}\,\text{elements}}\rangle
$$ with $r_j = r_j(t)\ge 0$.
Now we apply the formula $\langle A,B\rangle =\langle A\rangle
\langle B\rangle (1+[\overleftarrow{A}][\overrightarrow{B}]) $.
Note that in the case
$$
[\overleftarrow{A}]=[h_k,1,...],\,\,\,
[\overrightarrow{B}]=[1,h_{k+1},...]
$$
one has
$$
[\overleftarrow{A}] \ge\frac{1}{h_k+1}
\ge \frac{1}{h_{k+1}} ,\,\,\,
[\overrightarrow{B}]\ge \frac{1}{1+1/h_{k+1}}
=\frac{h_{k+1}}{h_{k+1}+1}
,\,\,\,
1+[\overleftarrow{A}][\overrightarrow{B}] \ge \frac{h_{k+1}+2}{h_{k+1}+1}.
$$
 So
 $$
 \langle 2, a_1, ..., a_t \rangle \ge
\prod_{k=1}^{f-1}\frac{h_{k+1}+2}{h_{k+1}+1}\times
 \langle
\underbrace{1,...,1}_{1+t-2(r_2+...+r_{n})}\rangle\langle
\underbrace{1,2,...,1,2}_{2r_2\,\text{elements}}\rangle
\cdots\langle
\underbrace{1,n,...,1,n}_{2r_{h_k}\,\text{elements}}\rangle
= $$
\begin{equation}
\prod_{k=1}^{f-1}\frac{h_{k+1}+2}{h_{k+1}+1}\times
 \langle
\underbrace{1,...,1}_{1+t-2(r_2+...+r_{n})}\rangle \times
\prod_{k=2}^{f-1} \langle
\underbrace{1,h_k,...,1,h_k}_{2r_k\,\text{elements}}\rangle .
 \label{L71}
\end{equation}
  One can easily see that
\begin{equation}
\langle \underbrace{1,...,1}_{l}\rangle \ge \lambda_1^{l-1}.
 \label{L72}
\end{equation}
Also one can see  from the continued fraction arguments that
\begin{equation}
\langle \underbrace{j,1,...,j,1}_{2r\,\text{elements}}\rangle= \frac{1}{\mu_j^{-1}-\mu_j} \left(
(1-\mu_j)\mu_j^r+(\mu_j^{-1}-1)\mu_j^{-r}\right) \ge \frac{j+\sqrt{j^2+4j}}{2\sqrt{j^2+4j}} \cdot \mu_j^r\ge \left( 1-\frac{1}{j}\right)
\cdot\mu_j^r
 \label{L73}
\end{equation}
where $\mu_j $ is defined in (\ref{muj}).
Moreover
$$
\prod_{k=1}^{f-1}\frac{h_{k+1}+2}{h_{k+1}+1}\times
\prod_{k=1}^{f-1} \frac{h_{k+1}-1}{h_{k+1}}
\ge \prod_{k=2}^\infty \left( 1-\frac{2}{k^2}\right)\ge \frac{1}{10}.
$$

Now (\ref{L7}) follows from (\ref{L71},\ref{L72},\ref{L73}). Lemma 12 is proved.

For   an irrational number $x=[a_1,a_2,...,a_t,...]$ for $ j=0,1,2,..$ define values $t_j$ inductively:
$$ t_0 =0,\,\,\, t_j=\min
\{t>t_{j-1}:\,\, a_t \ge 2\}.
$$
So $a_t \neq 1$ if and only if $ t = t_j$ for some $j\ge 1$. In other words all partial quotients between $ a_{t_j}$ and $a_{t_{j+1}}$ are equal
to one. Consider
$$
D_t = \max_{u\le t} (S_x(u)- \kappa_1 u).
$$
{\bf Lemma 13.}\,\,\,{\it Suppose $?'(x) =0$. Then for any $m$ large enough we have
\begin{equation}
t_{m+1}-t_m \le \frac{D_{t_m}}{\kappa_1-1}+1 ,
 \label{L81}
\end{equation}
 and for any function $\psi (t)$ under conditions
\begin{equation}
\psi (t) \ge 0,\,\,\,\, \psi (t) = o\left(\frac{\log\log t}{\log t}\right),\,\,\,\,
\psi (t)\cdot \log t \to +\infty,\,\,\,\,\,\, t\to\infty
\label{fff}
\end{equation}
one has
\begin{equation}
D_{t_{m}+1} \ge  \frac{ \sqrt{2\log \lambda_1-\log 2} }{\log 2}\cdot\sqrt{t_m\log t_m}\cdot (1
-t_m^{-\psi (t_m)})
.
 \label{L82}
\end{equation}
}

Proof.\,\,\,  Note that $ q_t(x) \ge \lambda_1^{t-1}$ for every $t$.  First of all as $?'(x) =0$ we deduce from  (\ref{12}) that there exists
$t_0(x) $ such that for each $t\ge t_0$  we have
$$
\frac{1}{\lambda_1^{3}} \ge
 \frac{\left|?\left( \frac{p_{t-1}}{q_{t-1}}\right) -?\left( \frac{p_{t}}{q_{t}}\right)\right|}{\left|
   \frac{p_{t-1}}{q_{t-1}}-\frac{p_{t}}{q_{t}}\right| }=
\frac{q_tq_{t-1}}{2^{S_x(t) }}  \ge  \frac{\lambda_1^{2t-3}}{2^{S_x(t)}}.
$$
  As    $
2^{\kappa_1}=\lambda_1^2$ we see that
\begin{equation}
S_x(t) -\kappa_1t \ge 0 \label{L85}
\end{equation}
for $t$ large enough. Particulary this inequality means that there
exist infinitely many  partial quotients $a_j$ greater than $1$.

Put $ c_{m+1} = t_{m+1} - t_m -1$.    From (\ref{12}) we have (by the same reasons) that the following inequality is valid for $m \ge m_0(x)$
(with some $m_0(x)$ depending on $x$ and large enough):
$$
\frac{1}{\lambda_1^{5}}\ge \frac{q_{t_{m+1}-2}(x)q_{t_{m+1}-1}(x)}{2^{S_x(t_{m+1}-1)}}\ge \frac{1}{\lambda_1^{5}}\cdot
\frac{\lambda_1^{2t_{m+1}}}{2^{\kappa_1t_m+D_{t_m}+c_{m+1}}}= \frac{1}{\lambda_1^{5}}\cdot\left(\frac{\lambda_1^2}{2}\right)^{c_{m+1}} \cdot
\frac{1}{2^{D_{t_m}}}.
$$
Now
$$
c_{m+1}\le \frac{D_{t_m}\log 2}{\log (\lambda_1^2/2)}=\frac{D_{t_m}}{\kappa_1-1}
$$
and  (\ref{L81}) follows.

By our notation the sequence $a_1,a_2,...,a_{t_m}$ is of the form
$$
a_1,a_2,...,a_{t_m}=
\underbrace{1,...,1}_{c_1},a_{t_1},\underbrace{1,...,1}_{c_2},a_{t_2},...
, \underbrace{1,...,1}_{c_m},a_{t_m},\,\,\,\,\,\, a_{t_j}\ge 2.
$$
 Now by (\ref{L85}) and by the definition of $D_t$ we have
\begin{equation}
a_{t_m} = S_x(t_m ) - S_x (t_m -1) = S_x(t_m )-\kappa_1 t_m - (S_x (t_m -1)-\kappa_1(t_m -1)) +\kappa_1 \le D_{t_m} +\kappa_1 \label{hvaa}
\end{equation}
 for
$m \ge m_0(x)$.

 As $ \kappa_1 \le 3/2$ we deduce that for $t$ large
enough we may suppose that $ r_1(t) >t/2$. So we can apply Lemma 12.
Remind that $
q_t\ge \frac{1}{8} \langle 2,a_1,...,a_t\rangle$.
With (\ref{12}) a for $m$ large enough it gives
$$
\frac{1}{80^2 \lambda_1^2}\ge \frac{q^2_{t_m}}{2^{S_x(t_m+1)}} \ge \frac{\lambda_1^{2t_m} \prod_{j=1}^m\left(\mu_{a_{t_j}}/\lambda_1^2\right)^2}{80^2 2^{S_{
x}(t_m+1)}}\ge
  \frac{\lambda_1^{2t_m} \prod_{j=1}^m\left(\mu_{a_{t_j}}/\lambda_1^2\right)^2}{80^2 2^{\kappa_1t_m+D_{t_m+1}+\kappa_1}}= \frac{
\prod_{j=1}^m\left(\mu_{a_{t_j}}/\lambda_1^2\right)^2}{80^2 2^{D_{t_m+1}+\kappa_1}}
 .
$$
So as $\lambda_1^2=2^{\kappa_1}$ we have
\begin{equation}
\sum_{1\le j\le m} \log (\mu_{a_{n_j}}/\lambda_1^2) \le \frac{\log 2}{2}\cdot  D_{t_m+1} . \label{hva1} \end{equation}

 As $a_{t_j} \ge 2$
we see that $\log (\mu_{a_{t_j}}/\lambda_1^2)\ge \log
(\mu_{2}/\lambda_1^2)> 1/3$.  So
$$
m \le \frac{3\log 2}{2}\cdot  D_{t_m+1}.$$

Now we take integer  $H$ to be large enough.
 From the definition (\ref{muj}) of $\mu_k$ we see that $ \mu_k \ge k $.
So
for $a_{t_j}\ge H$ for large $H$ we have
$$
\log (\mu_{a_{t_j}}/\lambda_1^2)= \log a_{t_j}-2\log\lambda_1+ \frac{\theta_H}{H},\,\,\,\,
 |\theta_H|\le 5.
$$
 So (\ref{hva1}) leads to
\begin{equation}
\sum_{1\le j\le m,\,\, a_{t_j}>H
 } \log a_{t_j} \le \frac{\log 2}{2}\cdot D_{t_m+1} +2m_1\log\lambda_1+\frac{\theta_H m_1}{H}
 \label{um2} \end{equation}
where
$$
m_1 = \#\{ j\in \{1,...,m\}:\,\,\, a_{t_j} \ge H\}.
$$
 Remind that
$$
  \sum_{j\le m}( a_{t_j} -\kappa_1) -(\kappa_1-1)\sum_{j\le m} c_j
=  S_x(t_m) - \kappa_1 t_m\ge 0 $$
by (\ref{L85}).
So
 $$t_m =
\sum_{j\le m} c_j + m
 \le \frac{1}{\kappa_1-1}\sum_{j\le m} (a_{t_j}-\kappa_1)+m\le
 \frac{1}{\kappa_1-1}\sum_{j\le m} a_{t_j}
.$$

Then from the definition of $m_1$ we have
\begin{equation}
t_m \le \frac{1}{\kappa_1-1} \sum_{j\le m} a_{t_j} \le \frac{1}{\kappa_1-1}\left(mH+m_1 \max_{1\le j\le
m} a_{t_j}\right). \label{ka}
\end{equation}
So for large $H$
 from (\ref{um2})  we have
$$
m_1 \log H \le \frac{\log 2}{2}\cdot D_{t_m+1} +2m_1\log\lambda_1+\frac{\theta_Hm_1}{H}.
$$
From the last inequality we see that for $H$ large enough the following estimate is valid:
$$
m_1\le  \frac{\log 2}{2}\cdot\frac{ D_{t_m+1}}{\log H} +
\frac{10 D_{t_m+1}}{(\log H)^2}.
$$

 From (\ref{hvaa}) and (\ref{ka}) we deduce that for large  $H $ we have
$$
t_{m} \le \frac{1}{\kappa_1-1} \left( \frac{3\log 2}{2} D_{t_m+1}H + m_1 \max\left( \max_{j\le m_0(x)}a_{t_j} \,\,;\,\,\kappa_1 + D_{t_m+1}\right)\right)
.$$
So $D_t\to\infty $ as $t\to\infty$ and
$$
t_m \le
\frac{\log 2}{2(\kappa_1-1)}\cdot\frac{D_{t_m+1}^2}{\log H}\cdot\left(
1+ \frac{3H\log H}{D_{t_m+1}}+\frac{20}{\log H}+\frac{20}{\log D_{t_m+1}}\right).
$$
  Inequality (\ref{L82}) follows from  the last inequality by taking
 $$
 H = D_{t_m+1}^{1-\psi (t_m)}/\log D_{t_m+1}
$$
(here we take into account conditions (\ref{fff})).
 Lemma is proved.

  We introduce some more notation. For $ r \in \{ 2,5,6,7,8,9,...\} $ put
\begin{equation}
\varphi_r(t) =\varphi_{r,x}(t)= \kappa_rt -S_x(t). \label{varphi}
\end{equation}

 {\bf Lemma 14.}\,\,\,{\it
1. Let
for an irrational $ x \in (0,1)$  for infinitely many values of  $t$ one has
\begin{equation}
4.1 t< S_x(t) < 4.9 t .\label{SX1}
\end{equation}
Then for these values of $t$ the following inequality is valid:
\begin{equation}
\langle 1,1,a_1,...,a_t \rangle =
 o\left( \left(
\frac{\lambda_4^{5}}{\lambda_5^4}\right)^{\frac{\varphi_2(t)}{\kappa_2}} 2^{\frac{S_x(t)}{2}}\right),\,\,\,\,\,
t \to \infty.
\label{notrn}
\end{equation}

2.
 Let
for some $t$ one has
\begin{equation}
4t \le S_x(t)\le5 t. \label{SX1p}
\end{equation}
Then
\begin{equation}
\langle 1,1,a_1,...,a_t \rangle \le
\left(\frac{\lambda_4^{5}}{\lambda_5^4}\right)^{\frac{\varphi_2(t)}{\kappa_2}} 2^{\frac{S_x(t)}{2}+3}.
\label{notrnp}
\end{equation}

}

Note that $ \frac{\lambda_4^5}{\lambda_5^4}>1$.

Proof of Lemma 14.

Let the inequality (\ref{SX1p}) is valid.
We shall use  the second statement of Lemma 5.
Put $s = S_x(t)$.
It follows that the maximal value of the continuant under consideration attains when its elements
are of the form $a$ and $a+1$. Obviously in our situation $a=4,a+1=5$. So
$$
\langle 1,1,a_1,...,a_t \rangle\le\langle
1,1,d_1,...,d_t\rangle ,\,\,\,\,\,
d_j \in \{ 4,5\},\,\,\,
d_1+...+d_t = a_1+...+a_t
.
$$
This inequality is not enough for our purpose. We should make use of the second statement of Lemma 6. It gives
$$
\langle 1,1,a_1,...,a_t \rangle\le \left(\frac{385}{384} \right)^{-\sigma_{2,t}/2} \langle 1,1,d_1,...,d_t\rangle
$$
with
$$
\sigma_{2,t} = \sigma_{2}(a_1,...,a_t) = \sum_{j=1}^t \min \{ |a_j -4|,|a_j -5|\}
 .$$
Applying Lemma 9 we get
$$
\langle 1,1,a_1,...,a_t \rangle\le   \left(1+\frac{1}{4\cdot 7^8}\right)^{{-\sigma_{2,t}/2}-\rho_t/2} \langle 1,1,
\underbrace{4,...,4}_{5t-S_x(t)},\underbrace{5,...,5}_{S_x(t)-4t}\rangle \le
$$
$$\le
8\cdot \left(1+\frac{1}{4\cdot 7^8}\right)^{{-\sigma_{2,t}/2}-\rho_t/2}
\langle\underbrace{4,...,4}_{5t-S_x(t)}\rangle\langle\underbrace{5,...,5}_{S_x(t)-4t}\rangle
$$
with
$$
\rho_t =\rho (d_1,...,d_t ) =
\sum_{j=1}^{t-1}|d_j-d_{j+1}|-2.
$$
Observe that
$$
\langle \underbrace{h,...,h}_{j}\rangle \le \lambda_h^j.
$$
Now
$$
\langle 1,1,a_1,...,a_t \rangle\le 8\cdot   \left(1+\frac{1}{4\cdot 7^8}\right)^{{-\sigma_{2,t}/2}-\rho_t/2}
\lambda_4^{5t-S_x(t)}\lambda_5^{S_x(t)-4t} .
$$
From the notation (\ref{varphi}) we see that
$$
\left( \frac{\lambda_5}{\lambda_4}\right)^{{S_x(t)}} \left( \frac{\lambda_4^{5}}{\lambda_5^4}\right)^{t} = \left(
\frac{\lambda_4^{5}}{\lambda_5^4}\right)^{\frac{\varphi_2(t)}{\kappa_2}} 2^{\frac{S_x(t)}{2}}.
$$
So
\begin{equation}
\langle 1,1,a_1,...,a_t \rangle\le  \left(1+\frac{1}{4\cdot 7^8}\right)^{{-\sigma_{2,t}/2}-\rho_t/2} \left(
\frac{\lambda_4^{5}}{\lambda_5^4}\right)^{\frac{\varphi_2(t)}{\kappa_2}} 2^{\frac{S_x(t)}{2}+3}
.
\label{L14}
\end{equation}
 As $\sigma_{2,t},\rho_t \ge 0$ the inequality (\ref{notrnp}) is proved.

If $\max \{ \sigma_{2,t},\rho_t\}\to\infty$ when $t\to \infty$   we have (\ref{notrn}).

Let now the inequality (\ref{SX1}) is valid.

Suppose that with some $M$ for infinitely many values of $t$ we have $\max \{ \sigma_{2,t},\rho_t\}\le M$. Then for such $t$  the number of
digits $a_j\neq 4, 5$ from the sequence $a_1,...,a_t$ is bounded by some constant $M_1$ depending only on $M$.
Moreover either the number of $a_j =4$ or the number of $a_j=5$ is bounded by a constant
$M_2$ depending only on $M$.
 It means that the continued fraction
expansion for $x$ is of the form  $ x = [b_1,...,b_r,\overline{4}] $ or $ x = [b_1,...,b_r,\overline{5}]. $ So $\lim_{ t\to \infty} \frac{S_x(t)}{t} =4\,\,\text{or}\,\, 5$.
It is impossible under the condition  (\ref{SX1}).

  Lemma 14 is proved.

{\bf Remark 1 to Lemma 14.}  Here we should note that the proof of
Lemma 14 relays on the algorithm for obtaining the sequence with
the maximal continuant from the sequence $1,1,a_1,...,a_t$. This
algorithm consists of two stages. The first one uses unit
variations only. It transforms the sequence $1,1,a_1,...,a_t$ into
the sequence $1,1,d_1,...,d_t$. (The second stage uses
substitutions only.) We may take a certain type of this procedure
to ensure that the sequence $1,1,d_1,...,d_t$ depends on the
initial sequence $1,1,a_1,...,a_t$ only. Moreover we take  each
unit variation  from the first stage to  transform a sequence
$1,1,...,a,...,b,...$  (or $1,1,...,b,...,a,...$) into the
sequence $1,1,...,a+1,...,b-1,...$  (or $1,1,...,b-1,...,a+1,...$)
where $a,b$ satisfies the conditions of Lemma 4. Hence in each
unit variation we have $\sigma_{2}(...,a+1,...,b-1,...)\le
\sigma_{2}(...,a,...,b,...)-1 $ (or
$\sigma_{2}(...,b-1,...,a+1,...)\le \sigma_{2}(...,b,...,a,...)-1
$). Now the algorithm of transforming the  sequence
$1,1,a_1,...,a_t$ into the sequence $1,1,d_1,...,d_t$ uses not
more than  $\sigma_2(a_1,...,a_t)$ unit variations.

{\bf Remark 2 to Lemma 14.} In fact in Lemma 14 we have proven that
\begin{equation}
\langle 1,1,a_1,...,a_t \rangle
\le \left(1+\frac{1}{4\cdot 7^8}\right)^{-\frac{\sigma_{2}(a_1,...,a_t)}{2}-\frac{\rho (d_1,...,d_t)}{2}} \left(
\frac{\lambda_4^{5}}{\lambda_5^4}\right)^{\frac{\varphi_2(t)}{\kappa_2}} 2^{\frac{S_x(t)}{2}+3},
\label{L14bis}
\end{equation}
where the sequence $d_1,...,d_t$ is obtained from the sequence $a_1,...,a_t$.

Now we must
investigate the  procedure of obtaining the minimal continuant from the continuant
$\langle 1,1, a_1,...,a_t\rangle$.
First of all we use unit variations to transform the sequence $1,1,a_1,...,a_t$ into the sequence
$1,1,d_1,...,d_t$ such that $ d_i \in \{1,n\}$
(as in Lemma 5, part 1 and lemma 6, part 1).
We take a certain type of such a procedure  to ensure that the sequence $d_1,...,d_t$  depends on the sequence
$a_1,...,a_t$ only.
Then we  apply substitutions from Lemma 11.

{\bf Lemma 15.}\,\,\, {\it \,\,\, Let $x=[a_1,...,a_t,...]\in E_n, n\ge 4$. Let
\begin{equation}
2S_x(t) \le (n+1)t. \label{SX}
\end{equation}
Then
the following statements are valid.

1.  If $ n \ge 5$
$$
\langle 1,1,a_1,...,a_t \rangle \ge C_n
\times\left(1+\frac{1}{4(n+2)^8}\right)^{ g_x(t)}\times
 \left(
\frac{\lambda_1^{n+1}}{\mu_n}\right)^{\frac{\varphi
_n(t)}{(n-1)\kappa_n}}\times 2^{\frac{S_x(t)}{2}}
$$ where
\begin{equation}
g_x(t) =\frac{\sigma_1(a_1,...,a_t)}{2}+\omega (d_1,...,d_t)
\label{geee}
\end{equation}
and
 $\sigma_1(a_1,...,a_t)$ ,  $\omega(d_1,...,d_t)$  are defined in (\ref{1w1}), (\ref{2w2}) respectively and
\begin{equation}
C_n= C_n' \times\left(1+\frac{1}{4(n+2)^8}\right)^{-\frac{n}{4}}\ge
\frac{\lambda_1^2}{24\mu_n}
 ,\,\,\,\,\,\,\,\,
C_n   ' =  \frac{1}{10}
\left(\frac{\mu_n}{\lambda_1^2}\right)^{\frac{2-n}{n-1}}\ge\frac{\lambda_1^2}{10\mu_n}.
\label{ceen}
\end{equation}

2. If $n=4$
then
\begin{equation}
\langle 1,1, a_1,...,a_t\rangle
\ge \frac{\lambda_1^2}{10}
\mu_4^{S_x(t)/5-1}.
\label{ceen1}
\end{equation}

}

 Proof.

First of all we apply Lemma 5 (statement 1).
Put $s = S_x(t)$.
 Then we can suppose that not more than one element $z=d_j$ differs from $ 1 ,n$, while values $t,
S_x(t)$ do not change. When we replace  element $z$ by $1$  the continuant $ \langle 1,1,d_1,..., z,...,d_t \rangle $ may decrease:
$$
\langle 1,1,d_1,..., z,...,d_t \rangle \ge \langle 1,1,d_1,...,
1,...,d_t \rangle.
$$
So we can suppose that {\it all} $d_j$ are from the set $1,n$ and
the value $t$ does not change while $S_x(t)$ should be replaced by
$S' = S_x(t) -z+1\ge S_x(t) -n+2   .$ Then for the values $w_j(t)$
we have
$$
w_2(t) = w_3(t) =...= w_n(t) = \frac{S'-t}{n-1} .
$$
Now we apply Corollary 3 of Lemma  12 (remember that $\mu_1=\lambda_1^2$ and that $\mu_j$ are increasing in $j$, so $\mu_n/\lambda_1^2 \ge 1$)
and see that
\begin{equation}
\langle 1,1,a_1,...,a_t \rangle \ge \frac{1}{10} \left(
\frac{\mu_n}{\lambda_1^2}\right)^{\frac{S'}{n-1}} \left(
\frac{\lambda_1^{n+1}}{\mu_n}\right)^{\frac{t}{n-1}} \ge \left(
\frac{1}{10} \left(\frac{\mu_n}{\lambda_1^2}\right)^{\frac{2-n}{n-1}}
\right)
  \left( \frac{\mu_n}{\lambda_1^2}\right)^{\frac{S_x(t)}{n-1}} \left( \frac{\lambda_1^{n+1}}{\mu_n}\right)^{\frac{t}{n-1}}
 .
\label{inequa}
\end{equation}
For $n=4$ we put (\ref{SX}) into (\ref{inequa}) and get (\ref{ceen1}).

Consider the case $n\ge 5$.
 By (\ref{varphi}) we have
$$
\left( \frac{\mu_n}{\lambda_1^2}\right)^{\frac{S_x(t)}{n-1}} \left( \frac{\lambda_1^{n+1}}{\mu_n}\right)^{\frac{t}{n-1}}= \left(
\frac{\lambda_1^{n+1}}{\mu_n}\right)^{\frac{\varphi _n(t)}{(n-1)\kappa_n}} 2^{\frac{S_x(t)}{2}}.$$
We have proven the estimate
 $$
\langle 1,1,a_1,...,a_t \rangle \ge C_n '\left(
\frac{\lambda_1^{n+1}}{\mu_n}\right)^{\frac{\varphi
_n(t)}{(n-1)\kappa_n}} 2^{\frac{S_x(t)}{2}}.
$$
Now we must  explain the arrival of an additional factor
\begin{equation}
\left(1+\frac{1}{4(n+2)^8}\right)^{\frac{2\sigma_1-n}{4}+\omega (d_1,...,d_t)}.
\label{factor}
\end{equation}
This factor appears from Lemmas 6 (part 1) and 11 as we can obtain the minimal continuant from a given one  applyind at first unit variation procedures (and these procedures lead to first additional factor from Lemma 6) and at second substitutions procedures (another additional factor arrives from Lemma 11). The factor (\ref{factor}) takes into account the infuence of the both two factors described behind.

Lemma is proved.

{\bf 11.\,\,\, Proofs of theorems.}

 {\bf 11.1.\,\,\, Proof of  Theorem 1 statement (i).}\,\,\
Let $x=[a_1,...,a_t,...]$. Note that $q_t(x) \ge 2^{-3} \langle
2,a_1,...,a_t\rangle $.
 By
   Lemma 1 it is sufficient to prove that under the conditions
 of  the statement (i) of Theorem 1 one has
 $$
 W_t(x):=
 \frac{\langle
2,a_1,...,a_t\rangle \langle 2,a_1,...,a_{t-1}\rangle}{{2}^{S_x(t)}} \to +\infty,\,\,\,\, t \to
 \infty.
 $$
Consider values $w_2(t)$ defined in (\ref{WE}). Note that
$w_2(t+1) \ge w_2(t)$ for any $t$. We distinguish two cases
\vskip+0.3cm
 {\bf Case 1.} \,\,\,$\kappa =
\limsup_{ t\to +\infty} \frac{S_x(t)}{t}
<\kappa_1
 $.
 \vskip+0.3cm
{\bf Case 2.}\,\,\, $\limsup_{ t\to +\infty} \frac{S_x(t)}{t}
=\kappa_1
$.
\vskip+0.3cm

 In the case 1
 by
Corollary 3 with $N=1$
 and the condition (\ref{eq1}) of Theorem 1 we see
that
$$
 W_t(x)\ge
\frac{1}{100\lambda_1\cdot 2^{C}} \times \left(\frac{\lambda_1^2}{{2^\kappa}}\right)^t\to  +\infty,\,\,\,\, t \to
 +\infty
 $$
(as $\lambda_1^2=2^{\kappa_1}>2^\kappa$).

 In the case 2 we see that $ w_2(t) \to +\infty$ as $t \to +\infty$. Also we should
take into account that under the conditions of Theorem 1 (i)  one has
for $t$
large enough
 $$
r_1(t) >\frac{t}{2},\,\,\, \sum_{i: 1\le i\le t,\, \, a_i\ge 2} a_i = S_x(t) -t+w_2.
$$
 Remind that the number of partial quotients  greater than 1 is $w_2 = r_2+...+r_n$ (here $r_j \ge 0$ and $n$ is the maximal partial quotient).
 Remind also that
 $ \mu_2<\mu_3<...<\mu_n$ and $\mu_j >j$.
Consider the value
$$
\mu = \min \mu_2^{r_2}\cdots \mu_n^{r_n}
$$
where the minimum is taken over the set
$$
\left \{ r_2,...,r_{n-1}\ge 0, r_n \ge 1:\,\, \sum_{j=2}^n r_j =
w_2,\,\, \sum_{j=2}^n jr_j = S_x(t)-t+w_2 \right\} .$$ We see that
$$
\mu \ge \mu_2^{w_2-1}\mu_n\ge
 \mu_2^{w_2-1}n\ge
  \mu_2^{w_2-1}\times
  \frac{\sum_{i: 1\le i\le t,\, \, a_i\ge 2} a_i}{w_2}
  \ge
 \mu_2^{w_2-1}\times
  \frac{S_x(t)-t+w_2}{w_2}
 .
  $$
Application of Lemma 12 gives
$$
\langle 2,a_1,...,a_t\rangle \ge
\frac{\lambda_1^{t}}{10}\prod_{k=2}^{n} \left(
\frac{\mu_k}{\lambda_1^2} \right)^{r_k(t)} \ge
\frac{\lambda_1^{t-2}}{10}\left(\frac{\mu_2}{\lambda_1^2}\right)^{w_2-1}
\frac{S_t(x) -t}{w_2}
$$
and of course
$$
 \langle 2,a_1,...,a_{t-1}\rangle \ge \lambda_1^{t}.
$$
So
$$
 W_t(x)\ge
\frac{\lambda_1^{2t}
\left(\frac{\mu_2}{\lambda_1^2}\right)^{w_2-1} (S_t(x)
-t)}{10\lambda_1^2w_2\cdot t\cdot2^{\kappa_1t+C}} \gg \frac{1}{w_2}\cdot
\left(\frac{\mu_2}{\lambda_1^2}\right)^{w_2}
 \to +\infty,\,\,\,\, t \to
 \infty
 $$
as $w_2\to +\infty$. Statement (i) of Theorem 1 is proved.

 {\bf 11.2.\,\,\, Proof of  Theorem 1 statement (ii).}\,\,\ We shall
 prove that for any function $\psi (t)$ increasing to infinity
 there exists an irrational number $x\in (0,1)$ such that
 $?'(x)$ does not exist  and
 $$
 S_x(t) \le \kappa_1 t+\frac{\log t}{\log 2 }+\psi (t).
 $$
We may suppose that $\psi (1) \ge 1 $ and that $\psi (t)=0(t),\,t \to \infty$.
We define integers $c_j,Q_j\, j=-1,0,1,2,...$ and $ t_j, b_j,
\,\, j=0,1 ,2,... $ by the following inductive procedure. Put
$c_{-1}=0,Q_{-1}=1, t_0=b_0=0.$
 Now  suppose that integers $c_j , q_j \,\, j=-1,0,1 ,2,..., k-1 $  and  $t_j, b_j,j\,\, j=0,1 ,2,..., k $ are defined. We must construct integers $c_k,Q_k,
t_{k+1}, b_{k+1}$. First of all take $c_k$ to be large enough to satisfy inequalities
\begin{equation}
\kappa_1 c_{k } >\sum_{j=0}^{k }c_j +\sum_{j=0}^{k }b_j
 +2\kappa_1, \label{c0}
\end{equation}
\begin{equation}
Q_{k-1}^2<2^{\frac{\psi (t_k+c_k+1)}{2}}
 , \label{c00}
\end{equation}
\begin{equation}
 c_k \ge \left[
\kappa_1(t_k+c_k+1)+\frac{\log (t_k+c_k+1)}{\log 2} +\psi (t_k+c+k+1)\right]-
\sum_{j=0}^kc_j-\sum_{j=0}^k b_j
. \label{c1}
\end{equation}
 It is possible to do this as in the right hand side of (\ref{c0}) the coefficient for $c_k$ is equal 1 meanwhile in the
 left hand side the corresponding coefficient is equal to $\kappa_1 >1$; as at the same time $Q_{k-1} $ from
 (\ref{c00}) does not depend on $c_k$. Moreover as $c_k$
increases the the right hand side of (\ref{c1}) depends on $c_k$
approximately as $(\kappa_1-1)c_k$.

 Then put
$$
t_{k+1} = t_{k}+c_{k}+1,\,\,\, b_{k+1} = \left[\kappa_1 t_{k+1}+\frac{\log t_{k+1}}{\log 2}+ \psi (t_{k+1})\right] -\sum_{j=1}^{k}c_j
-\sum_{j=1}^{k}b_j.
$$
From (\ref{c0}) we see that $ b_{k+1} \ge 2$. Now we define
 $$
Q_k =\langle  \underbrace{1,...,1}_{c_0}, b_1, \underbrace{1,...,1}_{c_1}, b_2,\underbrace{1,...,1}_{c_2}, b_3,..., \underbrace{1,...,1}_{c_{k-1}}, b_k,
\underbrace{1,...,1}_{c_{k}}, b_{k+1}\rangle.
$$
This continuant consist of $t_{k+1}$ partial quotients. Note that
\begin{equation}
 Q_{k-1}\lambda_1^{c_k-1} b_{k+1}\le
 Q_{k }
\le 4 Q_{k-1}\lambda_1^{c_k } b_{k+1}
 . \label{c3}
\end{equation}

 Now we take irrational  $x$ of the form
$$
x =[a_1,...,a_t,...]=[\underbrace{1,...,1}_{c_0}, b_1, \underbrace{1,...,1}_{c_1}, b_2, \underbrace{1,...,1}_{c_2}, b_3,..., \underbrace{1,...,1}_{c_{k-1}},
b_k,\underbrace{1,...,1}_{c_{k}}, b_{k+1},...] .$$
Now
$$
Q_k = q_{t_{k+1}}(x),\,\,\, k=1,2,3,...   .
$$
We see that for $t_k\le t<t_{k+1} $
$$
S_x(t)= \sum_{j=0}^{k-1} c_j + \sum_{j=1}^{k-1} b_j + b_k + t - t_k \le 1+\kappa_1 t_{k}+\frac{\log t_{k}}{\log 2}+ \psi (t_{k})+ t-t_k
 \le \kappa_1 t+\frac{\log t}{\log 2}+ \psi (t).
$$
So for $x$ constructed the inequality (\ref{eq2}) from the statement (ii) of Theorem 1 is true. At the same time from (\ref{12}) and (\ref{c1})
and the upper bound from (\ref{c3}) we see that
$$
 \frac{\left|?\left( \frac{p_{t_{k+1}-1}}{q_{t_{k+1}-1}}\right) -?\left( \frac{p_{t_{k+1}}}{q_{t_{k+1}}}\right)\right|}{\left|
   \frac{p_{t_{k+1}-1}}{q_{t_{k+1}-1}}-\frac{p_{t_{k+1}}}{q_{t_{k+1}}}\right| }=
\frac{q_{t_{k+1}}q_{t_{k+1}-1}}{2^{a_1+...+a_{t_{k+1}} }}\le
\frac{16Q_{k-1}^2\lambda_1^{2c_k}b_{k+1}}{2^{ \sum_{j=0}^{k} c_j +
\sum_{j=1}^{k+1} b_j }} \le
\frac{32Q_{k-1}^2\lambda_1^{2c_k}b_{k+1}}{2^{  \kappa_1t_{k+1}+\frac{\log
t_{k+1}}{\log 2}+\psi (t_{k+1})}} \le 2^{-\frac{\psi (t_{k+1})}{2}} \to 0
$$
as $ k \to
\infty .
$
From the other hand from (\ref{12}) we deduce by means of (\ref{c0}) and the lower bound from (\ref{c3}) the following inequality:
$$
  \frac{\left|?\left( \frac{p_{t_k-2}}{q_{t_k-2}}\right) -?\left( \frac{p_{t_{k}-1}}{q_{t_k-1}}\right)\right|}{\left|
   \frac{p_{t_k-2}}{q_{t_k-2}}-\frac{p_{t_k-1}}{q_{t_k-1}}\right| }
= \frac{q_{t_{k}-1}q_{t_{k}-2}}{2^{a_1+...+a_{t_{k}-1}}} \ge
$$
$$
   \ge
   \frac{(\prod_{j=1}^{k-1} b_j)\times\lambda_1^{2\sum_{j=0}^{k-1}c_j-2k}}{2^{ \sum_{j=0}^{k-1} c_j + \sum_{j=1}^{k-1} b_j }} \ge
   \frac{(\prod_{j=1}^{k-1} b_j)\times\lambda_1^{2\sum_{j=0}^{k-1}c_j-2(k-1)}}{2^{\kappa_1c_{k-1} }}
 \ge
   \prod_{j=1}^{k-1} b_j \to
\infty, \,\,\, k\to \infty .$$
So $?'(x)$ does not exist.

Statement (ii) of Theorem 1 is proved.

 {\bf 11.3.\,\,\, Proof of  Theorem 2 statement (i).}\,\,\
Let $t_m< t\le t_{m+1}$. Suppose $m$ to be large enough.
 We see from (\ref{L82}) of Lemma 13 that
\begin{equation}
D_t \ge
 D_{t_{m}+1} \ge
\frac{\sqrt{2\log \lambda_1 -\log 2}}{\log 2}\cdot
\sqrt{t_{m}\log t_{m}}\cdot (1-t_m^{-\psi (t_m)}).
\label{1111}
\end{equation}
We apply (\ref{L81})  of Lemma 13 to see that
\begin{equation}
 t_{m} \ge t - \frac{D_{t_m}}{\kappa_1-1} -1
\ge t - \frac{D_t}{\kappa_1-1} -1.
\label{2222}
\end{equation}
We substitute (\ref{2222}) into (\ref{1111}) and obtain the result of the statement (i) from Theorem 2.

 {\bf 11.4.\,\,\, Proof of  Theorem 2 statement (ii).}\,\,\,
 Put
$$  b_k =\left[\frac{2 k \log (k+3)}{\log 2}\right],\,\,\,\, k=1,2,3,...$$ and define
$$c_0= 2,\,\,\,c_{k} =\max \left (1; \left[
\frac{b_k -2\log (2^4b_k)/\log 2}{\kappa_1 -1} \right]\right) ,
$$
$$
t_k = \sum_{j=0}^{k} c_j + k.
$$
 Then as $\lambda_1^2/2 = 2^{\kappa_1-1}$ we have for all $ k$ large enough ($k \ge k_0 (x)$)
\begin{equation}
\left( \frac{\lambda_1^2}{2}\right)^{c_{k }} \le
\frac{2^{b_{k}}}{2^{8}b_{k}^2}.
 \label{gh}
\end{equation}
 We take    $x$ of the form
$$
x =[a_1,...,a_t,...]=[ \underbrace{1,...,1}_{c_0}, b_1,
\underbrace{1,...,1}_{c_1}, b_2, \underbrace{1,...,1}_{c_2},
b_3,..., \underbrace{1,...,1}_{c_{k-1}},
b_k,\underbrace{1,...,1}_{c_{k}}, b_{k+1},...] .$$ Now  as
$\lambda_1^2/2> 1$ for any $c\le c_k$ (with  $k \ge k_0(x)$) we see that
$$
\frac{q_{t_{k}-c}^2}{2^{S_x(t_{k}-c)}} \le
\frac{2^6q_{t_{k-1}}^2\lambda_1^{2(c_{k}-c)}b_k^2}{
2^{S_x(t_{k-1})+c_{k}+b_k-c}}\le \frac{2^6q_{t_{k-1}}^2 b_k^2}{
2^{S_x(t_{k-1}) }} \cdot \left(
\frac{\lambda_1^2}{2}\right)^{c_{k}}
 $$
and by (\ref{gh}) we have
$$
\frac{q_{t_{k}-c}^2}{2^{S_x(t_{k}-c)}} \le \frac{1}{2}\cdot
\frac{q_{t_{k-1}}^2}{2^{S_x(t_{k-1})}} \le \frac{1}{2^k}\cdot
\frac{q_{t_{k_0}}^2}{2^{S_x(t_{k_0})}}\to 0,\,\,\, k \to \infty.
$$
So
$$
\frac{ q_{r}^2}{2^{a_1+...+a_{r}  }} \to 0,\,\,\, r \to \infty
$$
 and we apply Lemma  2 to see that $?'(x) =0$.

Now we must prove an upper bound for
$S_x(t)$.
 By the definition of $c_j$ we have the following lower bound for the values of $t_k$:
$$
t_k  =\sum_{j=0}^{k} c_j + k
=\sum_{j=0}^{k_0-1}c_j +\sum_{j=k_0}^kc_j +k\ge
\frac{k^2\log k}{(\kappa_1-1)\log 2}\left(1-\frac{32}{k}\right).
$$
So
$$
k\le \sqrt{2(\kappa_1-1)\log 2}\cdot \sqrt{\frac{t_k}{\log t_k}}+ 64,\,\,\,
c_{k+1}
\le
\sqrt{\frac{2}{(\kappa_1-1)\log 2}}\cdot \sqrt{t_k\log t_k}+2^7\log t_k.
$$

From the definition of $c_j$ we deduce that
$$
c_{j} \ge \frac{1}{\kappa_1-1} \left( b_j-\frac{\log (16 b_j)}{\log\sqrt{2}} - 1\right)
$$
or
$$
b_j \le (\kappa_1-1)c_{j} + \frac{2}{\log 2} \log j + 4\log\log ( j+3).
$$
Now we combine all the estimates behind and deduce an upper bound
for $D_t$ from the range $ t_k < t\le t_{k+1}$:
$$
S_x(t)-\kappa_1 t
=S_x (t_{k+1}-(t_{k+1}-t)) -\kappa_1(t_{k+1}-(t_{k+1}-t))=
S_x(t_{k+1})-\kappa_1 t_{k+1} +(t_{k+1}-t)(\kappa_1-1)
\le$$
$$
\le S_x(t_{k+1})-\kappa_1 t_{k+1} +(\kappa_1-1)c_{k+1}
= \sum_{j=0}^{k}c_j +\sum_{j=1}^k b_j
-\kappa_1 \left(\sum_{j=0}^k c_j + k\right)
+(\kappa_1-1)c_{k+1}
 =$$
$$=
(1-\kappa_1) \sum_{j=0}^k c_j+ \sum_{j=1}^k b_j -\kappa_1 k +(\kappa_1-1)c_{k+1}
\le (\kappa_1-1)c_{k+1}+\frac{2k\log k}{\log 2} + 8k\log\log k  \le
$$
$$
\le
 2\sqrt{\frac{2(\kappa_1-1)}{\log 2}} \cdot \sqrt{t_k\log t_k} +2^6\sqrt{\frac{t_k}{\log t_k}}\log\log t_k\le
2\sqrt{\frac{2(\kappa_1-1)}{\log 2}} \cdot \sqrt{t\log t} +2^6\sqrt{\frac{t}{\log t}}\log\log t.
$$
Statement (ii) of
Theorem 2 is proved.

 {\bf 11.5.\,\,\, Proof of  Theorem 3 statement (i).}\,\,\,
By Lemma 2 it  is sufficient to prove that we have
\begin{equation}
q_t^2(x) = o(2^{S_x(t)}),\,\,\,\,\, t \to\infty.
\label{th3}
\end{equation}

{\bf Case 1.}\,\,\,  $S_x(t) < 4.9t$.

Remind that $\kappa_2>4.1$.

From (\ref{eq3}) it follows that
$\varphi_2 (t) \le C$. So the first statement of
 Lemma 14  leads to the bound (\ref{th3}).

{\bf Case 2.}\,\,\,  $4.9t \le S_x(t) \le 5t$.

Here $\varphi_2 (t) \le -0.4t$.
So the second statement of lemma 14 leads to the result of Theorem 3.

{\bf Case 3.}\,\,\,  $S_x(t) > 5t$.

Take $ t >100$.
 Consider a sequence
 \begin{equation}
 a_1,...,
a_{j-1}, a_j , a_{j+1}, ..., a_t \label{seqv}
\end{equation}
 with the sum
$$
a_1+...+a_j+...+a_t = S_x(t)
.$$
We may replace this sequence by the sequence
 $$
 a_1,...,
a_{j-1},1, a_j-1 , a_{j+1}, ..., a_t .$$ Then the sum of the digits does not change but the length of the sequence increases by $1$. Applying
this  replacement several times instead of the sequence (\ref{seqv}) we obtain a sequence
$$
b_1,....,b_{t'}
$$
with the same sum of digits
$$
S'(t') =  S_x(t)
$$
such that $$ S'(t') \le 5 t'.$$
As $t >100$ we see that $ S'(t')\ge 4.9 t'$.

Note that
$$
\langle  a_1,..., a_{j-1}, a_j , a_{j+1}, ..., a_t \rangle \le \langle  a_1,..., a_{j-1},1, a_j-1 , a_{j+1}, ..., a_t \rangle.
$$
So
$$
\langle  a_1, ..., a_t \rangle\le \langle  b_1,....,b_{t'}\rangle.
$$
Now we derive an upper bound for the continuant $ \langle 1,1, b_1,....,b_{t'}\rangle
$
from the second statement of
Lemma 14.
Now
  (\ref{th3}) follows analogously to the case 2.

Statement (i) of Theorem 3 is proved.

 {\bf 11.6.\,\,\, Proof of  Theorem 3 statement (ii).}\,\,\,
 We must construct a number $x$ such that $?'(x)$ does not exist but
(\ref{eq4}) holds  for all $t$. Put
\begin{equation}
x = [ \underbrace{5,...,5}_{c_1}, \underbrace{4,...,4}_{b_1}, \underbrace{5,...,5}_{c_2}, \underbrace{4,...,4}_{b_2},...,
\underbrace{5,...,5}_{c_k}, \underbrace{4,...,4}_{b_k},...
 ],
\label{54}
\end{equation}
\begin{equation}
t_k = \sum_{j=1}^k (c_j + b_j)
\label{544}
\end{equation}
where $c_k,b_k$ are defined by the following inductive procedure.

Let $t_0 = c_0=b_0 = 0$. Now suppose that $ c_0,b_0,...,c_{k-1}, b_{k-1}$ are defined. For a natural $c_k$ we define
$$
b_k= b (c_k) = \left[ \frac{(5-\kappa_2)c_k+\psi (t_{k-1}+c_k)-\psi (t_{k-1})}{\kappa_2-4} \right],
$$
$$
Q(t_k) =\langle  \underbrace{5,...,5}_{c_1}, \underbrace{4,...,4}_{b_1}, \underbrace{5,...,5}_{c_2}, \underbrace{4,...,4}_{b_2},...,
\underbrace{5,...,5}_{c_{k-1}}, \underbrace{4,...,4}_{b_{k-1}}, \underbrace{5,...,5}_{c_{k}}, \underbrace{4,...,4}_{b (c_k)} \rangle,
$$
 $$
S(t_k)= \sum_{j=1}^{k-1} (5c_j + 4b_j)+5c_k+ 4b(c_k).
$$

 We take $c_k$ to be large enough to satisfy
the condition
\begin{equation}
\frac{Q(t_k)^2}{2^{  S(t_k)}} \ge 1. \label{notex}
\end{equation}
  It is possible to do as from the definition (\ref{kappa2}) of $\kappa_2$ we have
$$
\frac{\lambda_5^{2(\kappa_2-4)}\lambda_4^{2(5-\kappa_2)}}{2^{\kappa_2}}=1
$$
and so
$$
\frac{Q(t_k)^2}{2^{  S(t_k)}} \gg \frac{\lambda_5^{2c_k}\lambda_4^{2b(c_k)}}{2^{5c_k +4b(c_k)}}\gg \left( \frac{\lambda_4^5}{\lambda_5^4}
\right)^{\frac{2\psi (t_{k-1}+c_k)}{\kappa_2}}\to \infty ,\,\,\, c_k\to \infty.
$$
 So given $ c_0,b_0,...,c_{k-1}, b_{k-1}$ we define $c_k, b_k = b(c_k)$. Real $x$ is defined by its continued fraction expansion.
So $q_{t_k}(x) = Q(t_k)$ and $ S_x(t_k) = S(t_k)$.

As for any $k$ we have (\ref{notex})   we see that
 the equality $?'(x) =0$ is not possible. It is not difficult to see from the construction that in the case
 $\psi (t) =o(\sqrt{t})$
the derivative $?'(x) $ {\it  does not} exist.
(In fact   it follows from Theorem 4.)

Consider $t$ from the interval
 $ t_{k-1}<t \le
t_k$. This interval can be divided into two intervals:
 $ t_{k-1}<t \le
t_{k-1}+c_k
$ (the first interval) and
 $ t_{k-1}+c_k<t \le
t_k=t_{k-1}+c_k+b_k$ (the second interval).  In the first interval as $ 5>\kappa_2$ and $\psi(\cdot)$ increases we have
$$
S_x(t) =S_x(t_{k-1})+5c\ge
\kappa_2t_{k-1}-\psi (t_{k-1})+5c\ge\kappa_2(t_{k-1}+c)-\psi (t_{k-1}+c)
$$
(here $c=t-t_{k-1}$)
and everything is fine. In the second interval for any $b\le b_k$ we have
$$
S_x(t_{k-1}+c_k+b) = S_x(t_{k-1})+5c_k+4b \ge
\kappa_2 t_{k-1}-\psi (t_{k-1}) +5c_k+4b=
$$
$$
=
\kappa_2(t_{k-1}+c_k+b) -\psi (t_{k-1})+ (5-\kappa_2)c_k-(\kappa_2-4)b\ge
\kappa_2(t_{k-1}+c_k+b)-\psi (t_{k-1}+c_k)\ge
$$
$$
\ge
\kappa_2(t_{k-1}+c_k+b)-\psi (t_{k-1}+c_k+b).
$$
Statement (ii) of Theorem 3 is proved.

{\bf 11.7.\,\,\, Proof of  Theorem 4 statement (i).}\,\,\,
As $?'(x)=+\infty$  we see by Lemma 2 that
$$
\frac{q_t(x)}{2^{\frac{S_x(t)}{2}}}\to \infty,\,\,\,
t\to \infty.
$$
So  by  formula (\ref{L14bis}) from Remark 2 to Lemma 14
$$
\varphi_2(t) \ge G
 (\sigma_2 (a_1,...,a_t)+\rho (d_1,...,d_t)),\,\,\,\,\,\,G=\frac{\kappa_2 \log (1+(4\cdot 7^8)^{-1})  }{2(5\log \lambda_4 - 4 \log \lambda_5)}
$$
for $t$ large enough
(here we use the notation from the proof   of Lemma 14, see also Remark 1 to Lemma 14).

If $\sigma_2 (a_1,...,a_t)+\rho (d_1,...,d_t)\ge \sqrt{t}/10$  then $
\varphi_2(t) \ge G \sqrt{t}/10$ and (\ref{theoche}) follows.

If $\sigma_2 (a_1,...,a_t)+\rho (d_1,...,d_t)< \sqrt{t}/10$ then
 $\rho (d_1,...,d_t)< \sqrt{t}/10$. Consider the sequence $1,1,d_1,...,d_t$ from the proof of Lemma 14.
We should note that $\kappa_2 > 4.4$ Hence the  number  of partial quotients among $d_1,...,d_t$ which are equal to $5$
is not less than $ 0.4t$. At the same time these partial quotients are distributed into not more than
$\rho (d_1,...,d_t) <\sqrt{t}/10$ blocks. So there exist $k \ge 0.4t /0.1\sqrt{t} = 4\sqrt{t}$ consecutive digits $d_{\nu+1},...,d_{\nu+k}$ equal to $5$.
So $\sigma_{2,t}<\sqrt{t}/10<k/40$.
 We consider the sequence
$$
1,1,d_1,...,d_\tau = 1,1,d_1,...,d_\nu ,\underbrace{5,...,5}_{k}
,\,\,\,\,\tau =\nu+k \le t.
$$
Remind that under the conditions of Theorem 4 we have $?'(x) =+\infty$.
From Lemma 2 for $t$ large enough we deduce the following inequalities (here we take into account Remark to the Definition of  the unit variation and Remark 1 to lemma 14):
$$
2^8\le
\frac{q^2_\tau (x)}{2^{S_x(\tau)}}\le
\frac{(2^{\sigma_2(a_1,...,a_t)}\times \langle 1,1,d_1,...,d_\nu ,\underbrace{5,...,5}_{k}\rangle )^2}{2^{S_x(\tau)}}\le
$$
$$
\le \frac{(2^{\sigma_2(a_1,...,a_t)}\times \langle 1,1,d_1,...,d_\nu ,\underbrace{5,...,5}_{k}\rangle )^2}{2^{d_1+...+d_\nu+5k-\sigma_2(a_1,...,a_t)}}\le
2^{3\sigma_2(a_1,...,a_t)+2}\cdot
\left(\frac{\lambda_5^2}{2^5}\right)^k\cdot
\frac{\langle 1,1,d_1,...,d_\nu \rangle^2}{2^{d_1+...+d_\nu}}.
$$
As $d_j \in \{4,5\}$  we have $4\nu \le d_1+...+d_\nu \le 5\nu$.
From the second part of
Lemma 14 we  have
$$
\frac{\langle 1,1,d_1,...,d_\nu \rangle^2}{2^{d_1+...+d_\nu}}\le2^6\cdot
\left(\frac{\lambda_4^5}{\lambda_5^4}\right)^{2\cdot \frac{\kappa_2\nu - d_1 -...-d_\nu}{\kappa_2}}.
$$
But
$$
\kappa_2\nu - (d_1+...+d_\nu) = \varphi_{2,x}(\nu) +\theta \sigma_2(a_1,...,a_t),\,\,\,\, |\theta |\le 1.
$$
Now
$$
2^8\le 2^{3\sigma_2(a_1,...,a_t)+8} \cdot
\left(\frac{\lambda_5^2}{2^5}\right)^k\cdot
\left(\frac{\lambda_4^5}{\lambda_5^4}\right)^{2(\varphi_{2,x}(\nu )+\sigma_2 (a_1,...,a_t))}
\le 2^{\frac{3k}{40}+8}\cdot\left(\frac{\lambda_5^2}{2^5}\right)^k
\cdot\left(\frac{\lambda_4^5}{\lambda_5^4}\right)^{2\varphi_2(\nu )+\frac{2k}{40}}
.$$
So
$$
1\le
\left(
\frac{\lambda_5^{72}\lambda_4^{10}}{2^{197}}\right)^k\cdot
\left(\frac{\lambda_4^{400}}{\lambda_5^{320}}\right)^{\varphi_2(\nu )}
$$
and $ \varphi_2 (\nu )\ge 0.0005 k \ge 0.002\sqrt{t}.$
So
Statement (i) of Theorem 4 is proved.

{\bf 11.8.\,\,\, Proof of  Theorem 4 statement (ii).}\,\,\,
We shall construct $x$ in the from (\ref{54}), $t_k$ should be defined as in (\ref{544}). Of course the choise of parameters $c_j,b_j$ will be different from that from the Theoerm
3 (statement (ii)).  Put
$$
c_1=100,\,\, b_k=100+10k ,\,\,\, c_{k+1}= \left[\frac{b_k(\log \lambda_4 -\log 4) - \log 60}{5\log\sqrt{2}-\log \lambda_5}\right],
$$
Then
$$
c_{k+1} = \frac{b_k(\log \lambda_4 -\log 4) }{5\log\sqrt{2}-\log \lambda_5}
-70\theta
 =\frac{\kappa_2-4}{5-\kappa_2}\times b_k -70 \theta ,\,\,\, 0\le \theta <1.
$$
$$
\left(\frac{\lambda_5^2}{2^5}\right)^{c_k-1}
\left(\frac{\lambda_4^2}{2^4}\right)^{b_{k-1}}>2,\,\,\,\,\,
\frac{\lambda_5^2}{2^5}<1<\frac{\lambda_4^2}{2^4}.
$$
Now
$$
\frac{
q^2_{t_{k-1}+c_k}(x)
}{2^{S_x(t_{k-1}+c_k)}}\ge
\frac{
q^2_{t_{k-2}+c_{k-1}}(x)
}{2^{S_x(t_{k-2}+c_{k-1})}}\times\left(\frac{\lambda_5^2}{2^5}\right)^{c_k-1}
\left(\frac{\lambda_4^2}{2^4}\right)^{b_{k-1}}
>2\frac{
q^2_{t_{k-2}+c_{k-1}}(x)
}{2^{S_x(t_{k-2}+c_{k-1})}}\gg 2^k\to +\infty
$$
as $ k \to \infty$.
But for $ c\le c_k$ we have
$$
\frac{
q^2_{t_{k-1}+c}(x)
}{2^{S_x(t_{k-1}+c)}}\ge
\frac{
q^2_{t_{k-1}+c_k}(x)
}{2^{S_x(t_{k-1}+c_k)}}
.
$$
Also for any $ b\ge 1$ we have
$$
\frac{
q^2_{t_{k-1}+c_k+b}(x)
}{2^{S_x(t_{k-1}+c_k+b)}}
\ge
\frac{
q^2_{t_{k-1}+c_k}(x)
}{2^{S_x(t_{k-1}+c_k)}}
\times
\left(\frac{\lambda_4^2}{2^4}\right)^{b}
\ge
\frac{
q^2_{t_{k-1}+c_k}(x)
}{2^{S_x(t_{k-1}+c_k)}}.
$$
So
$$
\frac{
q^2_{t}(x)
}{2^{S_x(t)}}\to+\infty,\,\,\, t\to \infty
$$
and $?'(x) =+\infty$.

Now we see that
$$
\varphi_2(t_k)
=\kappa_2t_k -S_x(t_k)\le \varphi_2(t_{k-1})+(\kappa_2-4)b_k-(5-\kappa_2)c_k=
$$
\begin{equation}
=
\varphi_2(t_{k-1})+(\kappa_2-4)b_{k-1}-(5-\kappa_2)c_k
+(\kappa_2-4)(b_k - b_{k-1})=\varphi_2(t_{k-1}) + 70 (5- \kappa_2 \theta )+ 10(\kappa_2-4)
 \label{nowwese}
\end{equation}
(here $\varphi_2(\cdot )$ is defined in (\ref{varphi})).

So we obtain recursive inequality
$$
\varphi_2(t_k)\le\varphi_2(t_{k-1})+400
$$
and $\varphi_2 (t_k) \le 400k$, while
$$
t_k \ge \sum_{j=1}^k b_j\ge 5k^2.
$$
From the other hand (\ref{nowwese}) leads to
$$
\varphi_2(t_{k}) >\varphi_2 (t_{k-1})
.
$$
Moreover in the example under consideration for any $k$  the function $d\mapsto \varphi _2(t_{k-1}+d)$
decreases in the interval $0\le d \le c_k$ and increases in the interval $c_k \le d \le c_k+b_k$.
It follows from the equality
$$
\varphi_2 (t_{k-1}+d) =\varphi_2 (t_{k-1}) - (5-\kappa_2)d,
\,\,\, 0\le d \le c_k,
$$
and from the equality
$$
\varphi_2 (t_{k-1}+d) =\varphi_2 (t_{k-1}+c_k)+ (\kappa_2 -4)(d-c_k),\,\,\,
c_k \le d \le c_k+b_k.
$$
Hence in the interval $ 0\le d\le c_k+b_k$
one has
$$
\varphi_2 (t_{k-1}+d) \le \varphi_2 (t_k).
$$

Let $t_{k-1}\le t<t_k$. then
$$
\varphi _2(t) \le \varphi_2(t_k)\le 400k\le 400\sqrt{t_k/5}\le 400 \sqrt{t /5}\cdot (t_k/t_{k-1}) \le 200\sqrt{t}
$$
for $t$ large enough. Theorem 4 is proved.

 {\bf 11.9.\,\,\, Proof of  Theorem 5 statement (i).}\,\,\,
Remind that $\kappa_n<\frac{n+1}{2}$ and  from (\ref{eq1n}) we have $\varphi_n(t) =\kappa_n t- S_x(t) \ge -C$. Note that the functions $g_x(t)$
defined in (\ref{geee})
 increases in $t$.

We consider two cases.

{\bf Case  1.} \,$g_x(t) \to \infty, \,\,\, t\to \infty$;

{\bf Case 2.} \,$ g_x(t) $ is bounded as $t\to \infty$.

In the case 1  the result follows from Lemmas 1 and 15 (part 1).

In the case 2 irrational $x$  has the following continued fraction
expansion: $x =[a_1,...,a_T,\overline{1}]$. So obviously $?'(x)
=+\infty$. Statement (i) of Theorem 5 is proved.

{\bf 11.10\,\,\, Proof of  Theorem 5 statement (ii).}\,\,\,
 The proof is close the the proof of Theorem 3 statement (ii).
It is necessary to
 construct a number $x$ such that $?'(x)$ does not exist but
(\ref{eq2n}) holds  for all $t$. Define
\begin{equation}
x = [ \underbrace{1,...,1}_{c_1}, \underbrace{1,n,1,n,...,1,n}_{2b_1 \,\,\text{digits}}, \underbrace{1,...,1}_{c_2}, \underbrace{1,n,1,n,...,1,n}_{2b_2\,\,\text{digits}},...,
\underbrace{1,...,1}_{c_k}, \underbrace{1,n,1,n,...,1,n}_{2b_k\,\,\text{digits}},...
 ],
\label{nadoobo}
\end{equation}
\begin{equation}
t_k = \sum_{j=1}^k (c_j +2 b_j).
\label{nadoobo1}
\end{equation}
Here
  $c_k,b_k$
we define inductively.
Put $t_0 = c_0=b_0 = 0$. Now suppose that $ c_0,b_0,...,c_{k-1}, b_{k-1}$ are defined. For a natural $c_k$ we define
$$
b_k= b (c_k) = \left[ \frac{(\kappa_n -1)c_k+\psi (t_{k-1}+c_k)-\psi (t_{k-1})}{n-2\kappa_n+1} \right].
$$
Now $$
Q(t_k) =
\langle
\underbrace{1,...,1}_{c_1}, \underbrace{1,n,1,n,...,1,n}_{2b_1 \,\,\text{digits}}, \underbrace{1,...,1}_{c_2}, \underbrace{1,n,1,n,...,1,n}_{2b_2\,\,\text{digits}},...,
\underbrace{1,...,1}_{c_k}, \underbrace{1,n,1,n,...,1,n}_{2b_k\,\,\text{digits}}
\rangle,
$$
  $$
S(t_k)= \sum_{j=1}^{k} (c_j + (n+1)b_j).
$$
Take $c_k$ to be large enough to satisfy
the condition
$$
\frac{Q(t_k)^2}{2^{  S(t_k)}} \le 1.
$$
  It is possible to do as from the definition (\ref{kappa1n}) of $\kappa_n$ we have
$$
\frac{\lambda_1^{2(n-2\kappa_n+1)}\mu_n^{\kappa_n-1}}{2^{(n-1)\kappa_n}}=1.
$$
 So one can easily see that in the case $\psi (t) = o(\sqrt{t})$
 the derivative $?'(x) $ does not exist (by Theorem 6).
 To prove that for all $t$ the inequality (\ref{eq2n}) is valid we need to perform the calculations similar to those from the proof of Theorem 3 statement (ii).

Statement (ii) of Theorem 5 is proved.

{\bf 11.11.\,\,\, Proof of  Theorem 6 statement (i).}\,\,\,

We shall give a sketched proof only. The proof follows the steps
of the proof of Theorem 4 statement (i) Suppose that all partial
qoutients of $x$ are bounded by $n$ and $?'(x)=0$. Analogously to
the function $\sigma_1(a_1,...,a_t)$ defined in Lemma 6  we must
consider a little bit more difficult function
$\sigma_1^{(n)}(a_1,...,a_t)$. The definition is as follows.

Given the sequence of partial qoutients $a_1,...,a_t$ we enumerate  all  of them different from $1, n$ in the non-decreasing order:
$$
1<
a_{i_1}\le a_{i_2}\le ...\le a_{i_k}<n,\,\,\,\, k\le t
$$
(so exactly $t-k$  partial quotients are equal to $1$ or $n$).
Put formally $a_{i_0}=1,a_{i_{k+1}}=n$.
Now to define $\sigma_1^{(n)}=\sigma_1^{(n)}(a_1,...,a_t)$ we put
$$
\sigma_1^{(n)}(a_1,...,a_t)
=
\sum_{j=0}^{k} (n-1-(a_{i_{j+1}}-a_{i_j})).
$$
We should note that $\sigma_1^{(n)}(a_1,...,a_t)\ge 0$ and in the case $k\ge 1$ we have
$\sigma_1^{(n)}(a_1,...,a_t)> 0$. Then analogously to (\ref{ggg})  Lemma 6 we deduce

\begin{equation}
\langle 1,1,a_1,a_2,...,a_t\rangle
 \ge \left( 1+\frac{1}{16(n+2)^3}\right)^{\frac{\sigma _1^{(n)}  -n+1}{2n-2}} \times
\min_{(b_1,b_2,...,b_t)\in U_n(t,s)} \langle
1,1,b_1,b_2,...,b_t\rangle . \label{ggg7}
\end{equation}

If  $S_x(t) \ge \frac{n+1}{2}t$ then for such $t$  theorem is
proven by $\frac{n+1}{2}>\kappa_n +\frac{1}{2} $. Moreover, by the
same reasons  theorem is valid in the case when there exists $\nu$
such that $\sqrt{t}\le \nu \le t$ and $S_x(\nu)\ge
\left(\frac{n+1}{2}-\frac{1}{10}\right)\nu $. It follows  from the
inequality
$$
\max_{u\le t} (S_x(u) -\kappa_n u)
\ge S_x(\nu ) -\kappa_n \nu \ge \left(\frac{n+1}{2} -\frac{1}{10}-\kappa_n\right) \nu
\ge \frac{\sqrt{t}}{3}.
$$
Hence we may assume that for all $\nu$  from the interval $\sqrt{t}\le \nu \le t$
we have
\begin{equation}
S_x(\nu ) < \left(\frac{n+1}{2} -\frac{1}{10}\right)\nu,\,\,\,\, \sqrt{t}\le \nu .
\label{buss}
\end{equation}
So the conditions of Lemmas 10, 11 are satisfied.
As
$
S_x(t) \le \frac{n+1}{2} t
$
and $ r_1(t) \ge \frac{t}{2}$
we may do the following.
Analogously to Lemma 15 we apply  (\ref{ggg7}) and Lemma 11 to get the inequality

\begin{equation}
\langle 1,1,a_1,...,a_t \rangle \ge C_n
\times\left(1+\frac{1}{4(n+2)^8}\right)^{ g^{(n)}_x(t)}\times
 \left(
\frac{\lambda_1^{n+1}}{\mu_n}\right)^{\frac{\varphi
_n(t)}{(n-1)\kappa_n}}\times 2^{\frac{S_x(t)}{2}}
\label{novee}
\end{equation}
  with
$$
g^{(n)}_x(t) =\frac{\sigma_1^{(n)}(a_1,...,a_t)}{2n-2}+\omega (d_1,...,d_t).
$$
Here $C_n$ is the same as in the formula (\ref{ceen})  from  Lemma
15. The algorithm from the proof of Lemma 6 (part 1) is also
modified: given a sequence $a_{i_1},...,a_{i_k}$ with  the
smallest element  $a_{i_1}>1$ and the  largest element $a_{i_k}<n$
we replace them by $a_{i_1}-1$ and $a_{i_k}+1$ correspondingly and
then enumerate all the elements of the new sequence (which are not
equal to $1$ or $n$) in non-decreasing order again. Moreover we
remark  here that at each step of the algorithm described the
value of $\sigma_1^{(n)}(a_1,...,a_t)$ decreases at least by 4 but
not more than $2n-2$. The sequence $d_1,...,d_t$ is just the
sequence appearing from $a_1,...,a_t$ after all unit variations.

Now form Lemma 1
and (\ref{novee})
we have for $t$ large enough
$$
-\varphi_n(t) \ge \frac{ g_x^{(n)}(t)}{16(n+2)^8}.
$$

Put
\begin{equation}
M=\frac{1}{10(n+2)^2}.
\label{mmee}
\end{equation}
Now we consider two cases.

If $g_x^{(n)}(t)\ge M\sqrt{t}$ then
$$
\max_{u\le t}(-\varphi_n(u))\ge
-\varphi_n(t) \ge
 \frac{ M}{16(n+2)^8}\times\sqrt{t}
$$
and theorem follows.

Consider the case $g_x^{(n)}(t)\le M\sqrt{t}$.

The proof of the formula (\ref{novee}) uses three stages. The
first step is the process of transformation of the initial
sequence of partial qoutients $a_1,...,a_t$ into the sequence
$d_1,...,d_t$ with elements $1$ and $n$ (with a possible exception
of one element) by means of a certain sequence of unit variations.
The second and the third stages are related to Lemma 11. The
second stage uses permutations to transform  the sequence
$d_1,...,d_t$ into the sequence in which  there is no consecutive
elements equal to $n$. The third stage collects together blocks of
the form $1,1,...,1$. In order to prove Theorem 6 we need to
consider what happens after the first and the second stages of the
process are completed.

 From the definition of $g_x^{(n)}(t)$ as $\omega (d_1,...,d_t) \ge -1$ we see that
$$
  \sigma_1^{(n)}(a_1,...,a_t)\le
(g^{(n)}_x(t) +1)(2n-2) \le
2C(n-1)\sqrt{t} + 2 n.
$$
So the first stage takes not more than
$\sigma_1^{(n)}/4 =
(M(n-1)\sqrt{t} +n)/2$ unit variations.

As
$$
\sum_{j\le t-1:\,\,d_j \ge 2}\delta (d_j,d_{j+1})\le
g_x^{(n)}(t)+1\le M\sqrt{t}+1
$$
we have done not more than
$M\sqrt{t}+1$
permutations during the second step of the process.

So we see that for  $t$ large enough
first two stages of the algorithm uses not more than
$(M(n-1)\sqrt{t}+n)/2+M\sqrt{t}+1\le Mn\sqrt{t}$
unit variations and substitutions.

Moreover,
$$
\sum_{j=1}^{t-2}
(1-\delta (d_j^*, d_{j+2}^*
)) \le
g_x^{(n)}(t)+1\le
 M\sqrt{t}+1.
$$
So
$$
\sum_{\sqrt{t}\le j\le t-2}
(1-\delta (d_j^*, d_{j+2}^*
)) \le M\sqrt{t}+1.
$$
Put $\gamma =[\sqrt{t}]+1$.
We see that after the first and the second stages are completed
the sequence of partial quotients
$a_\gamma,a_{
\gamma+1},..,a_t$
consist of not more than
$
M\sqrt{t}+2$ consecutive blocks of the form
$1,1,...,1
$ or $1,n,1,n,....,1,n$.
Note that there is not more  than $\frac{M\sqrt{t}+2+1}{2} \le M\sqrt{t}$ blocks $1,1,...,1$ among them (for $t$ large enough).

Remind that
$$
\kappa_n\le \kappa_5<2.5<\frac{n+1}{2}.
$$
So the total number of elements $1$ in all blocks of the form $1,1,...,1$
from the interval $[1,t]$
 is not less than
$t/10$. Hence the number of elements $1$ in blocks of the form
$1,...,1$ from the interval $[\gamma, t]$ is not less than
$t/10-\sqrt{t}$. These elements are located in not more than $
M\sqrt{t} $ blocks. Hence there exist a block $1,1,...,1$ of the
length
\begin{equation}
k \ge   \frac{ \sqrt{t}}{10M}  - \frac{1}{M}.
\label{odnam}
\end{equation}
Let now $d_1^*,...,d_t^*$ denotes the sequence of partial quotients which
appears after the first and the second stages of the process.
We consider the initial part of the sequence $d_1^*,...,d_t^*$  with  $ k$  units at the end:
$$
1,1,d_1^*,...,d_\tau^* \,\,=\,\, 1,1,d_1^*,...,d_\nu^*\underbrace{1,1,...,1}_k,\,\,\, \tau = \nu+k \le t,\,\,\,
\nu\ge \sqrt{t}.
$$
Note that every elementary procedure of the first and the second stage of the process (unit variation or permutation)
changes  the continuant $\langle 1,1,a_1,...,a_\nu\rangle$ or $\langle 1,1, d_1,...,d_\nu \rangle ,$
not more than by the factor $ 2(n+1)^2$. Really, for unit variations it follows from Remark after Definition 1;
for substitutions it follows from (\ref{knut}). To see this
one can take a sequence of the form
$$
C=
n,1,n,1,...,n,1
\,\,\,\,\,\,\text{or}\,\,\,\,\,\,
C= 1,n,1,n,...,1,n
$$
in the substitution $ K\to \Psi (A,B,C)$ from Lemma 7. Note that
this choice of the sequence $C$  leads to the following result.
Each unit variation or substitution during the first and the
second stages changes the value $S_x(\nu)$ not more than by $n-1$.
The total number of procedures does not exceed $Mn\sqrt{t}$. Hence
$|S_x(\nu) - (d_1^*+...+d_\nu ^*)|\le Mn^2\sqrt{t}$. So $\kappa_n
\nu - (d_1^*+...+d_\nu^*) = \varphi_{n,x}(\nu )+\theta \cdot
Mn^2\sqrt{t},\,\,|\theta |\le 1$. The conditions of Lemma 15 (part
1) are still valid for the number $x^*=[d_1^*,...,d_\nu^*,...]$
which appears from $x=[a_1,...,a_\nu ,...]$ after the first and
the second stages are completed. This fact follows from
(\ref{buss})  and the inequality $ Mn^2 \le 1/10$ as
$$
d_1^*+...+d_\nu^* < \left(\frac{n+1}{2}-\frac{1}{10}\right)\nu+Mn^2\sqrt{t}\le\left(
\frac{n+1}{2}-\frac{1}{10}+Mn^2\right)\nu \le \frac{n+1}{2}\nu
.
$$

As
$?'(x)=0$ Lemma 1
and the above arguments
lead to
$$
\frac{1}{16\lambda_1^4
\cdot (\frac{\lambda_1^2}{2})^{1/M}}
\ge
\frac{(q_{\tau-1})^2}{2^{S_x(\tau )}}\ge
\frac{\langle 1,1,a_1,..., a_{\tau -1}\rangle^2}{2^{S_x(\tau ) +6}}
\ge
\frac{\langle 1,1,d_1^*,...,d_\nu^*,\underbrace{1,1,...,1}_{k-1}\rangle^2}{
(2(n+1)^2)^{Mn\sqrt{t}}\times
2^{6+d_1^*+...+d_\nu^* +k)+Mn^2\sqrt{t}}}
\ge
$$
$$
\ge
\frac{4}{2^6\lambda_1^4}\times
\frac{1}{n^{4Mn\sqrt{t}}\cdot 2^{Mn^2\sqrt{t}}}\times
\left(\frac{\lambda_1^2}{2}\right)^k \times \frac{\langle d_1^*,...,d_\nu^*\rangle^2}{2^{d_1^*+...+d_\nu^* }}
.$$
We apply (\ref{odnam}) and Lemma 15 (part 1) for the sequence $d_1^*,...,d_\nu^*  $ and obtain
the inequality
$$
1\ge
\frac{(\lambda_1^2/2)^{\frac{\sqrt{t}}{10M}}}{4^{Mn^2\sqrt{t}}}\times
\left(
\frac{\lambda_1^{n+1}}{\mu_n}\right)^{2\frac{\kappa_n\nu - (d_1^*+...+d_\nu^*)}{(n-1)\kappa_n}}
\ge \frac{(\lambda_1^2/2)^{\sqrt{t}/(10M)}
}{4^{Mn^2\sqrt{t}}}\times
\left(\frac{\lambda_1^{n+1}}{\mu_n}\right)^{2\frac{\varphi_{n,x}(\nu ) - Mn^2\sqrt{t}}{(n-1)\kappa_n}}
.
$$
The last inequality leads to
$$
-\varphi_{n,x}(\nu ) \ge \frac{\sqrt{t}}{(n+2)^3}.
$$
Theorem 6 statement (i) is proved.

{\bf 11.12.\,\,\, Proof of  Theorem 6 statement (ii).}\,\,\

We shall construct $x$ in the form (\ref{nadoobo}). Values of
$t_k$ should be defined as in (\ref{nadoobo1}). Of course the
choice of parameters $c_j,b_j$ will be different of them from the
Theorem 5 statement (ii).

Put
$$
c_1=100,\,\,\,
b_k = 100+10k,\,\,\,
c_{k+1} =
\left[
\frac{b_k((n+1)\log 2-2\log \mu_n)-10}{2\log \lambda_1-\log 2}\right].
$$
Then
$$
c_{k+1} =
\frac{b_k((n+1)\log 2-2\log \mu_n)}{2\log \lambda_1-\log 2}
-70\theta
\le\frac{n-2\kappa_n+1}{\kappa_n -1}b_k - 70 \theta
,\,\,\,
0\le \theta \le 1.
$$
So
\begin{equation}
4\left(\frac{\lambda_1}{\sqrt{2}}\right)^{c_k}\left(\frac{\mu_n}{\sqrt{2}^{n+1}}\right)^{b_{k-1}}\le\frac{1}{2},\,\,\,
\frac{\mu_n}{\sqrt{2}^{n+1}}<1<\frac{\lambda_1}{\sqrt{2}}.
\label{maii}
\end{equation}
In the sequel the proof follows the steps of the proof of Theorem 4 statement (ii).
Inequalities (\ref{maii})  lead to
$\frac{q_t^2(x)}{2^{S_x(t)}}\to 0, t\to \infty$. So $?'(x) = 0$.

Moreover, for $\varphi_n(t)$ we have
$$
-\varphi_n(t_k) = S_x(t_k) -\kappa_n (t_k) =
-\varphi_n (t_{k-1}) + (n+1-2\kappa_n)b_k -(\kappa_n-1)c_k  =
$$
$$
=
-\varphi_n(t_{k-1})+(n+1-2\kappa_n)b_{k-1}
-(\kappa_n-1)c_k + (n+1-2\kappa_n)(b_k - b_{k-1})\le -\varphi_n (t_{k-1}) +30n.
$$
So we have recursive inequality
$$
-\varphi_n(t_k) \le -\varphi_n (t_{k-1}) + 30 n$$
and
$$
-\varphi_n(t_k)\le 30 kn.
$$
At the same time
$$
t_k \ge \sum_{j=1}^kb_j \ge 5k^2.
$$
The monotonicity of the function $d\mapsto -\varphi_n (t_{k-1}+d)$
in the interval $ 0\le d\le c_k + 2b_k$ is proved by the same arguments as in the proof of Theorem 4 statement (ii).
As a result we have
$$
-\varphi_n(t_{k-1}+d) \le - \varphi_n (t_k).
$$
Hence
 for $ t_{k-1} \le t < t_k$ we have
$$
-\varphi_n (t) \le -\varphi_n (t_k) \le 30 kn \le
30 n \sqrt{t_k/5}
\le
30n \sqrt{t/5} \,t_k/t_{k-1}
\le 15 n\sqrt{t}
$$
for $t$ large enough. Theorem 6 statement (ii) is proved.

{\bf 11.13.\,\,\, Proof of  Theorem 7.}\,\,\,

It is sufficient to deduce a lower bound for the continued fraction denominator
\begin{equation}
q^2_t(x) \ge (2+\varepsilon)^{S_x(t)}
\label{tys}
\end{equation}
for some positive $\varepsilon$.
By Lemma 5 statement 1 we may assume that all partial quotients are equal $1 $  or $4$. Applying Lemma 7 we see that
$$
q_t(x) \ge \frac{1}{8}\cdot
\langle 1,1, d_1,...,d_k,\underbrace{4,...,4}_{m},b_1,...,b_l\rangle
$$
for certain $m = m(t), k = k(t), l=l(t)$
where sequences $d_1,...,d_k$ and $ b_1,...,b_l$ consist of elements $1$ and $4$; moreover
there is no two consecutive $4$ in both of these sequences. Then
$$
\langle \underbrace{4,...,4}_{m}\rangle \ge \lambda_4^{m-1}
.
$$
For continuants
$$
\langle d_1,...,d_k\rangle,\,\,\,
\langle b_1,...,b_l\rangle
$$
the inequality  (\ref{SX}) is true and
we apply formula (\ref{inequa}) from Lemma 15 (part 2).
So
$$
\langle d_1,...,d_k\rangle \gg \mu_4^{\frac{d_1+...+d_k}{5}},\,\,\,
\langle  b_1,...,b_l\rangle\gg \mu_4^{\frac{b_1+...+b_l}{5}}
.$$
Now
$$
q_t(x)
\gg
\mu_4^{\frac{d_1+...+d_k+b_1+..+b_l}{5}}\lambda_4^{m}
.$$
As $\mu_4^{1/5} >1.42>\sqrt{2}, \lambda_4  >4.2 $
we have (\ref{tys}) with $\varepsilon = 0.001$. Theorem 7 is proved.


\begin{thebibliography}{99}

\bibitem{MI}
Minkowski H.\, Gesammelte Abhandlungen vol.2 (1911).
\bibitem{SA}

Salem R.\,\,\, On some singular monotone functions which are  strictly increasing. // Trans. Amer. Math. Soc., 53 (1943), 427 - 439.

\bibitem{KI}
 Kinney J.R. \,
Note on a singular function of Minkowski. // Proc. Amer. Math. Soc. 11 (1960), p. 788 - 789.


\bibitem{PARA}
 Paradis J.,  Viader P., Bibiloni L.\,
 A new light on Minkowski's $?(x)$ function. // J. Number Theory., 73 (1998), 212 -227.

\bibitem{PARA2}
 Paradis J.,  Viader P., Bibiloni L.\,
 The derivative of Minkowski's $?(x)$ function. // J. Math. Anal. and Appl. 253 (2001), 107 - 125.

\bibitem{arxiv}
Anna A. Dushistova, Nikolai G. Moshchevitin, On the derivative of the Minkowski question mark function $?(x)$.  // Preprint, available at arXiv:
0706.2219v2 [math.NT] 17Dec2007

\bibitem{knuth}
Graham R.L.,  Knuth D.E., Patashnik O.\,
Concrete Mathematics. Addison-Wesley, 1994.


\bibitem{knuth1}
 Knuth D.E. \,
Art of Computer Programming V.2., Third edition. Addison-Wesley, 1997.
\bibitem{KAN}

Kan I.D,\, Refinding of the comparison rule for continuants. //  Discrete Math. Appl. 10  (2000), no. 5, 477 - 480.


\bibitem{MOZ}

Motzkin T.S., Straus E.G. \, Some combinatorial extremum problems. // Proc. Amer. Math. Soc. (1956), 7, 1014 - 1021.




          \end{thebibliography}
\end{document}